\documentclass[12pt,leqno]{article}

\def\re{\mathop{\rm Re}}
\def\im{\mathop{\rm Im}}

\def\diam{\mathop{\rm diam}}
\def\span{\mathop{\rm span}}

\newtheorem{theorem}{Theorem}
\newtheorem{lemma}[theorem]{Lemma}
\newtheorem{proposition}[theorem]{Proposition}
\newtheorem{definition}[theorem]{Definition}
\newtheorem{corollary}[theorem]{Corollary}

\newcommand{\begintheorem}{\addtocounter{equation}{1}\begin{theorem}}
\newcommand{\beginlemma}{\addtocounter{equation}{1}\begin{lemma}}
\newcommand{\beginproposition}{\addtocounter{equation}{1}\begin{proposition}}
\newcommand{\begindefinition}{\addtocounter{equation}{1}\begin{definition}}
\newcommand{\begincorollary}{\addtocounter{equation}{1}\begin{corollary}}

\renewcommand{\thetheorem}{\arabic{section}.\arabic{equation}}
\renewcommand{\theequation}{\arabic{section}.\arabic{equation}}

\begin{document}

\title{Potpourri}

\author{Stephen William Semmes	\\
	Rice University		\\
	Houston, Texas}

\date{}

\maketitle


\renewcommand{\thefootnote}{}   

\footnotetext{These notes are connected to the ``potpourri'' course in
the mathematics department at Rice University in the fall of 2004, and
in particular they are very much influenced by the participants and
the discussions in class.}

\tableofcontents

\bigskip

	Let us begin with a mathematical personality quiz.  Consider
the five classes of functions given by polynomials, power series that
are entire or converge everywhere, power series which converge on the
unit disk, power series which converge on a disk of some positive
radius, and formal power series.  Each of these classes contains the
previous one.  For the quiz one should arrange them according to
preference.

\section{Real and complex numbers}
\label{real and complex numbers}
\setcounter{equation}{0}

	Let ${\bf Q}$, ${\bf R}$, and ${\bf C}$ denote the fields of
rational, real, and complex numbers, respectively.  As usual, a
complex number $z$ can be expressed as $x + y \, i$, where $x$ and $y$
are real numbers, and $i^2 = -1$, and we call $x$, $y$ the real and
imaginary parts of $z$ and denote them $\re z$, $\im z$, respectively.
We write ${\bf Z}$ for the integers, and ${\bf Z}_+$ for the set of
positive integers.

	Let us recall that the rational numbers are dense in the real
numbers in the sense that if $x$, $y$ are real numbers such that
$x < y$, then there is a rational number $r$ such that
\begin{equation}
	x < r < y.
\end{equation}
Also, for each positive real number $x$ there is a positive integer
$n$ such that
\begin{equation}
	n \, x \ge 1,
\end{equation}
which is the same as saying that
\begin{equation}
	\frac{1}{n} < x.
\end{equation}

	If $z = x + y \, i$ is a complex number, $x, y \in {\bf R}$,
then the \emph{complex conjugate} of $z$ is denoted $\overline{z}$ and
defined to be $x - y \, i$.
Thus
\begin{equation}
	2 \re z = z + \overline{z}
\end{equation}
and
\begin{equation}
	2 \, i \im z = z - \overline{z}
\end{equation}
for all $z \in {\bf C}$.  Furthermore,
\begin{equation}
	\overline{z + w} = \overline{z} + \overline{w}
\end{equation}
and
\begin{equation}
	\overline{z \, w} = \overline{z} \, \overline{w}
\end{equation}
for any complex numbers $z$, $w$.

	If $x$ is a real number, then the \emph{absolute value} of $x$
is denoted $|x|$ and is defined to be equal to $x$ when $x \ge 0$ and
to be equal to $-x$ when $x \le 0$.  Thus $|x|$ is always a
nonnegative real number which is equal to $0$ if and only if $x = 0$.
The triangle inequality for the absolute value function states that
\begin{equation}
	|x + y| \le |x| + |y|
\end{equation}
for all $x, y \in {\bf R}$, and this is easy to verify from the
definitions.  We also have that
\begin{equation}
	|x \, y| = |x| \, |y|
\end{equation}
for all $x, y \in {\bf R}$.

	If $z = x + y \, i$ is a complex number, with $x, y \in {\bf
R}$, then the norm or \emph{modulus} of $z$ is denoted $|z|$ and
defined by
\begin{equation}
	|z| = \sqrt{x^2 + y^2}.
\end{equation}
If $z$ happens to be a real number, then this is the same as the
absolute value of $z$ as in the preceding paragraph.  We can also
describe $|z|$ as the nonnegative real number such that
\begin{equation}
	|z|^2 = z \, \overline{z}.
\end{equation}
Notice that $|z| = 0$ if and only if $z = 0$.  Also,
\begin{equation}
	|\re z|, |\im z| \le |z|
\end{equation}
for any complex number $z$.

	If $z$, $w$ are complex numbers, then
\begin{equation}
	|z + w|^2 = (z + w) (\overline{z} + \overline{w})
	  	= z \, \overline{z} + z \, \overline{w} 
			+ \overline{z} \, w + w \, \overline{w}
\end{equation}
and therefore
\begin{eqnarray}
	|z + w|^2  & = & |z|^2 + 2 \re (z \, \overline{w}) + |w|^2	\\
	  & \le & |z|^2 + 2 |z| \, |w| + |w|^2 = (|z| + |w|)^2.
							\nonumber
\end{eqnarray}
In other words,
\begin{equation}
	|z + w| \le |z| + |w|
\end{equation}
for all complex numbers $z$, $w$.  We also have that
\begin{equation}
	|z \, w| = |z| \, |w|
\end{equation}
for all $z, w \in {\bf C}$, since the complex conjugate of a product
is equal to the product of the corresponding complex conjugates.

	The real numbers form a metric space with the standard metric
$|x - y|$, and the complex numbers form a metric space with the
standard metric $|z - w|$.  One can view the real numbers as a
subspace of the complex numbers, since the standard metric on ${\bf
R}$ is the same as the restriction of the standard metric on ${\bf C}$
to ${\bf R}$.  Notice that the set of rational numbers is dense inside
the real line in the sense of metric spaces, and that the complex
numbers with rational real and imaginary parts are dense in the
complex plane.

	If $\{z_j\}_{j=1}^\infty$, $\{w_j\}_{j=1}^\infty$ are
sequences of complex numbers which converge to $z, w \in {\bf C}$,
then
\begin{equation}
	\lim_{j \to \infty} z_j + w_j = z + w
\end{equation}
and
\begin{equation}
	\lim_{j \to \infty} z_j \, w_j = z \, w.
\end{equation}
If $\{z_j\}_{j=1}^\infty$ is a sequence of nonzero complex
numbers which converges to $z \in {\bf C}$, $z \ne 0$, then
\begin{equation}
	\lim_{j \to \infty} \frac{1}{z_j} = \frac{1}{z}.
\end{equation}
For any sequence of complex numbers $\{z_j\}_{j=1}^\infty$
which converges to a complex number $z$, we have that
\begin{equation}
	\lim_{j \to \infty} |z_j| = |z|,
\end{equation}
as one can show using the fact that
\begin{equation}
	\Bigl| |a| - |b| \Bigr| \le |a - b|
\end{equation}
for all complex numbers $a$, $b$, which is a consequence of the
triangle inequality.  Of course the same statements hold for the real
numbers as a special case.

	Let $\{z_j\}_{j=1}^\infty$ be a sequence of complex numbers.
One can check that $\{z_j\}_{j=1}^\infty$ is a Cauchy sequence if and
only if the sequences
\begin{equation}
	\{\re z_j\}_{j=1}^\infty, \quad \{\im z_j\}_{j=1}^\infty
\end{equation}
of real and imaginary parts of the $z_j$'s are Cauchy sequences as
sequences of real numbers.  If $z$ is a complex number, then
$\{z_j\}_{j=1}^\infty$ converges to $z$ if and only if $\{\re
z_j\}_{j=1}^\infty$, $\{\im z_j \}_{j=1}^\infty$ converge to $\re z$,
$\im z$, respectively, as sequences of real numbers.  A basic result
states that the real line with the standard metric is complete as a
metric space, which is to say that every Cauchy sequence of real
numbers converges.  It follows that the complex numbers are also
complete as a metric space.

	There is another notion of completeness for the real numbers,
which is based on ordering.  Suppose that $A$ is a subset of the real
line.  A real number $b$ is said to be an \emph{upper bound} for $A$
is $a \le b$ for all $a \in A$.  A real number $c$ is said to be a
\emph{least upper bound} or \emph{supremum} of $A$ if $c$ is an upper
bound of $A$ and if $c \le b$ for every upper bound $b$ of $A$.  It is
easy to see that the supremum of $A$ is unique if it exists.

	As an ordered set, the real numbers are complete in the sense
that every nonempty set $A$ of real numbers which has an upper bound
has a least upper bound.  One often starts with this and derives
completeness in the sense of convergence of Cauchy sequences.  One can
also start with completeness in the sense of convergence of Cauchy
sequences and derive completeness in the sense of ordering.

	A sequence $\{x_j\}_{j=1}^\infty$ of real numbers is said to
be \emph{monotone increasing} if $x_j \le x_{j+1}$ for all positive
integers $j$.  If there is an upper bound for the $x_j$'s, then the
sequence $\{x_j\}_{j=1}^\infty$ converges, and the limit is the
supremum of the $x_j$'s.  If one starts with completeness of the real
numbers in terms of ordering, then this statement follows easily from
the definitions, and at any rate the convergence of monotone
increasing sequences of real numbers which are bounded from above is
basically an equivalent form of completeness.

	Notice that every bounded subset of the real or complex
numbers is totally bounded, which is to say that it can be expressed
as the union of finitely many subsets of arbitrarily small diameter.
Basically this reduces to the fact that every interval in the real
line can be expressed as the union of finitely many subintervals of
arbtrarily small length.  It follows from the completeness of the real
and complex numbers as metric spaces that a subset of the real or
complex numbers is compact if and only if it is closed and bounded.

	Let $\sum_{j=0}^\infty a_j$ be an infinite series with terms
$a_j \in {\bf C}$.  We say that this series \emph{converges} if the
sequence of partial sums $s_l = \sum_{j=0}^l a_j$ converges, in which
case $\sum_{j=0}^\infty a_j$ is defined to be $\lim_{l \to \infty}
s_l$.  Of course the sequence of partial sums converges if and only if
it is a Cauchy sequence, which is equivalent to saying that
$\sum_{j=0}^\infty a_j$ converges if and only if for each $\epsilon >
0$ there is an $L \ge 0$ so that $\biggl | \sum_{j=l}^m a_j \biggr| <
\epsilon$ whenever $m \ge l \ge L$.  This is known as the \emph{Cauchy
criterion} for convergence of $\sum_{j=0}^\infty a_j$, and it follows
that if $\sum_{j=0}^\infty a_j$ converges, then $\lim_{j \to \infty}
a_j = 0$.  If $\sum_{j=0}^\infty a_j$ is an infinite series of
nonnegative real numbers, then the sequence of partial sums is a
monotone increasing sequence of nonnegative real numbers, and hence
converges if and only if the sequence of partial sums has an upper
bound.
	
	An infinite series $\sum_{j=0}^\infty$ of complex numbers is
said to \emph{converge absolutely} if the series $\sum_{j=0}^\infty
|a_j|$ converges.  Absolute convergence implies ordinary convergence,
because of the Cauchy criterion.  A series $\sum_{j=0}^\infty a_j$ of
complex numbers converges absolutely if and only if $\sum_{j=0}^\infty
\theta_j \, a_j$ converges for any sequence $\theta_0, \theta_1,
\ldots$ of complex numbers such that $|\theta_j| \le 1$ for all $j$.
If $\sum_{j=0}^\infty a_j$ and $\sum_{j=0}^\infty b_j$ are convergent
series of complex numbers and $\alpha$, $\beta$ are complex numbers,
then $\sum_{j=0}^\infty (\alpha \, a_j + \beta \, b_j)$ also
converges, and if $\sum_{j=0}^\infty a_j$, $\sum_{j=0}^\infty b_j$
converge absolutely, then $\sum_{j=0}^\infty (\alpha \, a_j + \beta \,
b_j)$ converges absolutely too.  A series which converges but does not
converge absolutely is said to converge conditionally.

	Let $z$ be a complex number, and consider the series
$\sum_{j=0}^\infty z^j$, where as usual $z^j$ is interpreted as being
equal to $1$ when $j = 0$, even when $z = 0$.  If $|z| \ge 1$, then
$|z|^j \ge 1$ for all $j$.  When $|z| < 1$, it is well-known that
$\lim_{j \to \infty} z^j = 0$.  For each positive integer $n$ we have
that $\sum_{j=0}^n z^j$ is equal to $(1 - z^{n+1}) / (1 - z)$ assuming
$z \ne 1$.  It follows that $\sum_{j=0}^\infty z^j$ converges when
$|z| < 1$, and in fact converges absolutely, with the sum equal to $1/
(1 - z)$.

	The \emph{Leibniz alternating series test} states that if
$b_0, b_1, b_2, \ldots$ is a monotone decreasing sequence of
nonnegative real numbers which converges to $0$, then the series
$\sum_{j=0}^\infty (-1)^j \, b_j$ converges.  This can be verified
using the Cauchy criterion.  More generally, suppose that the $b_j$'s
have the same property and that we have a sequence of complex numbers
$c_0, c_1, \ldots$ such that the partial sums $\sum_{j=0}^n c_j$ are
bounded.  In this case one can again show that $\sum_{j=0}^\infty c_j
\, b_j$ converges.  These results can be used to give examples of
series which converge but do not converge absolutely.

	If $a_0, a_1, \ldots$, are complex numbers, then we get the
associated power series $\sum_{j =0}^\infty a_j \, z^j$.  If this
series converges for some particular $z_0 \in {\bf C}$, then $\lim_{j
\to \infty} a_j \, z_0^j = 0$, and the sequence of $a_j \, z_0^j$'s is
bounded.  In this case one can check that $\sum_{j=0}^\infty a_j \,
z^j$ converges absolutely for all complex numbers $z$ such that $|z| <
|z_0|$, and $\sum_{j=0}^\infty a_j \, z^j$ converges absolutely when
$|z| \le |z_0|$ if $\sum_{j=0}^\infty a_j \, z_0^j$ converges
absolutely.  The \emph{radius of convergence} $R$ of
$\sum_{j=0}^\infty a_j \, z^j$, $0 \le R \le \infty$, is characterized
by saying that $\sum_{j=0}^\infty a_j \, z_j$ converges absolutely
when $|z| < R$ and does not converge at all when $|z| > R$.  Any power
series converges at $0$, and for a complex number $z$ with $|z|$ equal
to the radius of convergence of the series it may be that the series
does not converge, or converges conditionally, or converges
absolutely.

\section{$p$-Adic numbers}
\label{p-adic numbers}
\setcounter{equation}{0}

	Fix a prime number $p$, which is to say a positive
integer $p$ which is divisible only by $1$ and itself.
The \emph{$p$-adic absolute value} of a rational number $x$
is denoted $|x|_p$ and defined as follows.  If $x = 0$, then
$|x|_p = 0$.  If $x \ne 0$, then $x$ can be expressed as
\begin{equation}
	p^k \frac{m}{n}
\end{equation}
for some integer $k$ and nonzero integers $m$, $n$, where
neither $m$ nor $n$ is an integer multiple of $p$, and
one puts
\begin{equation}
	|x|_p = p^{-k}.
\end{equation}
Thus $|x|_p$ is small if $x$ has a lot of factors of $p$,
and it is large if $x$ has a lot of factors of $1/p$.

	Thus $|x|_p$ is always a nonnegative real number and $|x|_p =
0$ if and only if $x = 0$.  The $p$-adic absolute value satisfies a
stronger version of the triangle inequality, called the
\emph{ultrametric} version, which is that
\begin{equation}
	|x + y|_p \le \max (|x|_p, |y|_p)
\end{equation}
for all $x, y \in {\bf Q}$.  This is not too difficult to verify, and
we also have that
\begin{equation}
	|x \, y |_p = |x|_p \, |y|_p
\end{equation}
for all $x, y \in {\bf Q}$.

	Fix a prime number $p$.  The field of \emph{$p$-adic numbers}
is denoted ${\bf Q}_p$.  Basically ${\bf Q}_p$ completes the rational
numbers with respect to the $p$-adic metric in the same way that the
real numbers complete the rational numbers with respect to the
standard metric.  To be more precise, the $p$-adic numbers ${\bf Q}_p$
contain a copy of the rational numbers ${\bf Q}$.  The $p$-adic
absolute value function extends to ${\bf Q}_p$, is also denoted
$|\cdot |_p$, and satisfies the properties that $|x|_p$ is equal to
$0$ when $x = 0$ and is equal to $p^l$ for some integer $l$ when $x
\ne 0$.

	As on ${\bf Q}$, we have that $|x + y| \le \max (|x|_p,
|y|_p)$ and $|x \, y|_p = |x|_p \, |y|_p$ for all $x, y \in {\bf
Q}_p$.  The function $|x - y|_p$ defines an ultrametric on ${\bf
Q}_p$, extending the $p$-adic metric on ${\bf Q}$.  With respect to
this ultrametric, ${\bf Q}$ is a dense subset of ${\bf Q}_p$ and ${\bf
Q}_p$ is complete.

	Let $\{x_j\}_{j=1}^\infty$ and $\{y_j\}_{j=1}^\infty$ be
sequences of $p$-adic numbers which converge to $x, y \in {\bf Q}_p$,
respectively.  In this event we have that
\begin{equation}
	\lim_{j \to \infty} x_j + y_j = x + y
\end{equation}
and
\begin{equation}
	\lim_{j \to \infty} x_j \, y_j = x \, y
\end{equation}
just as for real and complex numbers.  If $x_j \ne 0$ for all $j$ and
$x \ne 0$, then
\begin{equation}
	\lim_{j \to \infty} x_j^{-1} = x^{-1}.
\end{equation}

	Let $\sum_{j=0}^\infty a_j$ be an infinite series whose terms
are $p$-adic numbers.  Just as for series of real or complex numbers,
we say that the series converges if the sequence of partial sums
converges in ${\bf Q}_p$.  If $\sum_{j=0}^\infty a_j$ is an infinite
series of $p$-adic numbers which converges and $\alpha$ is a $p$-adic
number, then $\sum_{j=0}^\infty \alpha \, a_j$ also converges.  If
$\sum_{j=0}^\infty a_j$, $\sum_{j=0}^\infty b_j$ are infinite series
of $p$-adic numbers which converge, then $\sum_{j=0}^\infty (a_j +
b_j)$ converges too.

	Because ${\bf Q}_p$ is complete, an infinite series
$\sum_{j=0}^\infty a_j$ of $p$-adic numbers converges if and only if
the sequence of partial sums forms a Cauchy sequence.  In the $p$-adic
case this is equivalent to
\begin{equation}
	\lim_{j \to \infty} a_j = 0,
\end{equation}
because of the ultrametric property of the $p$-adic absolute value
function.  In particular, in the $p$-adic case, if an infinite
series $\sum_{j=0}^\infty a_j$ converges, then so does
$\sum_{j=0}^\infty \theta_j \, a_j$ whenever $|\theta_j|_p \le 1$
for all $j$.

	Let $a_0, a_1, \ldots$, be a sequence of $p$-adic numbers, and
consider the associated power series $\sum_{j=0}^\infty a_j \, x^j$,
where again $x^j$ is interpreted as being equal to $1$ when $j = 0$
for all $x$.  If $\sum_{j=0}^\infty a_j \, x_0^j$ converges for some
particular $x_0 \in {\bf Q}_p$, then
\begin{equation}
	\lim_{j \to \infty} a_j \, x_0^j = 0
\end{equation}
in ${\bf Q}_p$, which is to say that
\begin{equation}
	\lim_{j \to \infty} |a_j|_p \, |x_0|_p^j = 0
\end{equation}
as a limit of real numbers.  It follows that
\begin{equation}
	\lim_{j \to \infty} a_j \, x^j = 0
\end{equation}
for all $p$-adic numbers $x$ such that $|x|_p \le |x_0|_p$.

	Let $x$ be a $p$-adic number.  If $x \ne 1$, then for each
positive integer $n$ we have that $\sum_{j=0}^n x^j$ is equal to $(1 -
x^{j+1}) / (1 - x)$.  If $|x|_p < 1$, then we get that
$\sum_{j=0}^\infty x^j$ converges, and that the sum is equal to $1 /
(1 - x)$.

	By definition of the $p$-adic absolute value, if $x \in {\bf
Z}$, then $|x|_p \le 1$.  Let $x$ be a rational number such that
$|x|_p \le 1$.  This means that we can write $x$ as $a/n$ for some
positive integer $n$ which is not a multiple of $p$.  More precisely
we can write $x$ as $a / (b + p \, l)$, where $b$, $l$ are integers
and $1 \le b < p$.  We may as well assume that $b = 1$, because
otherwise we could multiply both the numerator and denominator by a
positive integer $c$ such that $b \, c$ is equivalent to $1$ modulo
$p$.

	Thus we have a rational number $x$ which is expressed as $a /
(1 + p \, l)$ for some integers $a$, $l$.  It follows from the earlier
discussion that
\begin{equation}
	x = a \sum_{j=0}^\infty (-p \, l)^j,
\end{equation}
where the series converges in the $p$-adic metric.  Using the
approximation by partial sums we get that $x$ is a limit of integers
in the $p$-adic metric, so that the set of rational numbers with
$p$-adic absolute value less than or equal to $1$ is the same as the
closure of the set of integers in the $p$-adic metric.
Similarly the set of $p$-adic numbers with $p$-adic absolute
value less than or equal to $1$ is equal to the closure
of the set of integers as a subset of ${\bf Q}_p$ with respect
to the $p$-adic metric, and this set is denoted ${\bf Z}_p$
and called the $p$-adic integers.

	For each positive integer $n$, if $x$ is an integer, then we
can write $x$ as $b + p^n \, y$, where $b$, $y$ are integers and $0
\le b < p^n$.  In other words, not only does every integer have
$p$-adic absolute value less than or equal to $1$, but we can express
${\bf Z}$ as the union of $p^n$ subsets each with $p$-adic diameter
equal to $p^{-n}$.  In fact the same is true of the $p$-adic integer
${\bf Z}_p$, by approximation.  It follows that ${\bf Z}$, ${\bf Z}_p$
are totally bounded with respect to the $p$-adic metric.  Because
${\bf Q}_p$ is complete with respect to the $p$-adic metric and ${\bf
Z}_p$ is closed, we obtain that ${\bf Z}_p$ is a compact subset of
${\bf Q}_p$ with respect to the $p$-adic metric.

	Similarly, $p^{-l} \, {\bf Z}_p$ is a compact subset of ${\bf
Q}_p$ for each positive integer $l$.  To be more precise, $p^{-l} \,
{\bf Z}_p$ consists of the $x \in {\bf Q}_p$ of the form $p^{-l} \, y$
for some $y \in {\bf Z}_p$, which is the same as the set of $x \in
{\bf Q}_p$ such that $|x|_p \le p^l$.  Every bounded subset of ${\bf
Q}_p$ is contained in $p^{-l} \, {\bf Z}_p$ for some positive integer
$l$, and therefore a subset of ${\bf Q}_p$ is compact if and only if
it is closed and bounded.

\section{${\bf Z}[1/p]$}
\label{Z[1/p]}
\setcounter{equation}{0}

	Let ${\bf Z}[1/p]$ denote the set of rational numbers of the
form $p^{-l} \, a$, where $a$ is an integer and $l$ is a nonnegative
integer.  Clearly ${\bf Z}[1/p]$ is dense as a subset of ${\bf Q}$
with respect to the standard metric.  One can also check that ${\bf
Z}[1/p]$ is dense as a subset of ${\bf Q}$ with respect to the
$p$-adic metric.  This reduces to the fact that every rational number
$x$ with $|x|_p \le 1$ can be approximated by integers in the $p$-adic
metric.

	Consider the Cartesian product ${\bf Q} \times {\bf Q}$,
consisting of ordered pairs $(x_1, x_2)$ with $x_1, x_2 \in {\bf Q}$.
Let us combine the standard and $p$-adic metrics on ${\bf Q}$ into a
product metric on ${\bf Q} \times {\bf Q}$, in which the distance from
$(x_1, x_2)$ to $(y_1, y_2)$ in ${\bf Q} \times {\bf Q}$ is defined to
be the maximum of $|x_1 - y_1|$ and $|x_2 - y_2|_p$.  That is, we use
the standard distance in the first coordinate and the $p$-adic
distance in the second coordinate.

	The diagonal embedding of ${\bf Q}$ into ${\bf Q} \times {\bf
Q}$ sends $x \in {\bf Q}$ to $(x, x)$.  We can use this embedding to
map ${\bf Z}[1/p]$ into ${\bf Q} \times {\bf Q}$.

	Suppose that $x \in {\bf Z}[1/p]$.  If $|x|_p \le 1$, then $x
\in {\bf Z}$, in which case either $x = 0$ or $|x| \ge 1$.
In other words, either $x = 0$, or
\begin{equation}
	\max(|x|, |x|_p) \ge 1.
\end{equation}
If $x, y \in {\bf Z}[1/p]$, then either $x = y$, or
\begin{equation}
	\max(|x - y|, |x - y|_p) \ge 1.
\end{equation}
Thus the image of ${\bf Z}[1/p]$ in ${\bf Q} \times {\bf Q}$ is
discrete with respect to the product metric.  Namely, the distance
between any two distinct points in the image of ${\bf Z}[1/p]$ in
${\bf Q} \times {\bf Q}$ is at least $1$.

	Moreover, any point in ${\bf Q} \times {\bf Q}$ is at bounded
distance from a point in the image of ${\bf Z}[1/p]$ under the
diagonal embedding.  Explictly, if $(y, w)$ is any element of ${\bf Q}
\times {\bf Q}$, then there is an $x \in {\bf Z}[1/p]$ such that
\begin{equation}
	|x - y| < 1
\end{equation}
and
\begin{equation}
	|x - w|_p \le 1.
\end{equation}
We may as well assume that $y, w \in {\bf Z}[1/p]$, because ${\bf
Z}[1/p]$ is dense in ${\bf Q}$ with respect to both the standard and
$p$-adic metrics.  Let us write $w$ as $y + a + b$, where $a$ is an
integer and $0 \le b < 1$, and put $x = y + b$.  Then $x - y = b$ and
$x - w = -a$ have the required properties.

	Now let $E$ be a finite set of primes, which we can also
enumerate as $p_1, \ldots, p_n$, and let ${\bf Z}_E$ denote the set of
rational numbers of the form
\begin{equation}
	p_1^{-l_1} \cdots p_n^{-l_n} \, a,
\end{equation}
where $a$ is an integer and $l_1, \ldots, l_n$ are nonnegative
integers.  Thus ${\bf Z}_E$ is dense as a subset of ${\bf Q}$ with
respect to the standard metric $|x - y|$ as well as the $p_i$-adic
metrics $|x - y|_{p_i}$ for $i = 1, \ldots, n$.

	Consider the Cartesian product
\begin{equation}
	{\bf Q} \times {\bf Q} \times \cdots \times {\bf Q},
\end{equation}
with $n + 1$ copies of ${\bf Q}$.  We define the distance between two
points in this Cartesian product to be the maximum of the standard
distance between the first coordinates and the $p_i$-adic distance
between the $(i + 1)$th coordinate when $1 \le i \le n$.  We can embed
${\bf Z}_E$ into this Cartesian product using the diagonal embedding,
which sends $x \in {\bf Q}$ to an $(n+1)$-tuple whose coordinates are
all equal to $x$.

	If $x \in {\bf Z}_E$ and the $p_i$-adic absolute value of $x$
is less than or equal to $1$ for $i = 1, \ldots, n$, then $x$ is an
integer.  As a result, either $x = 0$, or the standard absolute value
of $x$ is greater than or equal to $1$.  This implies that if we take
two distinct elements of ${\bf Z}_E$ and consider their embeddings
into the Cartesian product of $n + 1$ copies of ${\bf Q}$, then the
distance between the two points in the Cartesian product is greater
than or equal to $1$.  Thus ${\bf Z}_E$ becomes discrete in the
Cartesian product.

	Suppose that $(y, w_1, \ldots, w_n)$ is an element of the
Cartesian product of $n + 1$ copies of ${\bf Q}$.  We would like to
show that there is a point in the image of ${\bf Z}_E$ under the
diagonal embedding whose distance to $(y, w_1, \ldots, w_n)$ is less
than or equal to $n$.  Specifically, let us check that there is an $x
\in {\bf Z}_E$ such that
\begin{equation}
	|x - y| < n
\end{equation}
and
\begin{equation}
	|x - w_i|_{p_i} \le 1
\end{equation}
for each $i$, $1 \le i \le n$.

	We may as well assume that $y \in {\bf Z}_E$, since ${\bf
Z}_E$ is dense in ${\bf Q}$ with respect to the standard metric.  For
$i = 1, \ldots, n$ we may assume that $w_i - y \in {\bf Z}[1/p_i]$,
since ${\bf Z}[1/p_i]$ is dense in ${\bf Q}$ with respect to the
$p_i$-adic metric.  Thus we can write $w_i$ as $y + a_i + b_i$, where
$a_i$ is an integer and $b_i \in {\bf Z}[1/p_i]$ satisfies $0 \le b_i
< 1$.  If we take $x = y + b_1 + \cdots + b_n$, then it is easy to see
that $x$ has the required properties.  This uses the fact that every
element of ${\bf Z}[1/q]$ has $p$-adic absolute value less than or
equal to $1$ when $p$, $q$ are distinct prime numbers.

\section{Exponential functions}
\label{exponential functions}
\setcounter{equation}{0}

	Let $\sum_{j=0}^\infty a_j \, z^j$ and $\sum_{l=0}^\infty b_l
\, z^l$ be formal power series.  If we multiply these two series
formally, then we get a power series $\sum_{n=0}^\infty c_n \, z^n$,
where
\begin{equation}
	c_n = \sum_{j=0}^n a_j \, b_{n-j}
\end{equation}
for each $n \ge 0$.  The sequence of $c_n$'s is called the
\emph{Cauchy product} of the $a_j$'s and $b_l$'s.

	Now suppose that $\sum_{j=0}^\infty a_j$ and
$\sum_{l=0}^\infty b_l$ are convergent series of complex numbers.  We
can define the $c_n$'s as in the preceding paragraph, and consider the
series $\sum_{n=0}^\infty c_n$.  Does this series necessarily
converge?  If so, is the sum equal to the product of the sums of the
$a_j$'s and $b_l$'s?

	If there are only finitely many nonzero $a_j$'s and $b_l$'s,
then this is simply an exercise in arithmetic.  Suppose now that the
$a_j$'s and $b_l$'s are nonnegative real numbers.  It is easy to see
that each partial sum of $\sum_{n=0}^\infty c_n$ is less than or equal
to
\begin{equation}
	\biggl(\sum_{j=0}^\infty a_j \biggr)
		\biggl(\sum_{l=0}^\infty b_l \biggr).
\end{equation}
Hence $\sum_{n=0}^\infty c_n$ converges, and the sum is less than or
equal to the aforementioned product.  One can show too that
$\sum_{n=0}^\infty c_n$ is equal to the product of the sums of the
$a_j$'s and $b_l$'s, because it is greater than or equal to the
product of any of their partial sums.

	Using this one can check that if $\sum_{j=0}^\infty a_j$ and
$\sum_{l=0}^\infty b_l$ converge absolutely, then $\sum_{n=0}^\infty
c_n$ converges absolutely.  Namely, one applies the previous case to
$|a_j|$, $|b_l|$, and one notes that $|c_n|$ is less than or equal to
the $n$th term of the Cauchy product of the absolute values of the
$a_j$'s and $b_l$'s.  One way to show that the sum of the $c_n$'s is
equal to the product of the sums of the $a_j$'s and the $b_l$'s is to
decompose the series into linear combinations of series with
nonnegative entries and apply the result already known for them.
Another way is to approximate the series of $a_j$'s and $b_l$'s by
finite sums.  For finite sums we get the right answer by arithmetic,
and the point is to check that small errors for the sums of the
$a_j$'s and $b_l$'s lead to small errors for the sum of $c_n$'s in a
suitable manner.

	It is a nice exercise to check that $\sum_{n=0}^\infty c_n$
converges, and that the sum is equal to the product of the sums of the
$a_j$'s and $b_l$'s, if one of $\sum_{j=0}^\infty a_j$ and
$\sum_{l=0}^\infty b_l$ has only finitely many terms and the other is
a convergent series.  A refinement of this states that if one of
$\sum_{j=0}^\infty a_j$, $\sum_{l=0}^\infty b_l$ converges absolutely
and the other converges, then $\sum_{n=0}^\infty c_n$ converges and is
equal to the product of the sums of the $a_j$'s and $b_l$'s.  See
\cite{Rudin2}.

	A theorem of Abel states that if $\sum_{j=0}^\infty a_j$,
$\sum_{l=0}^\infty b_l$, and $\sum_{n=0}^\infty c_n$ all converge,
then the sum of the $c_n$'s is equal to the product of the sums of the
$a_j$'s and $b_l$'s.  To prove this, let $r$ be a positive real number
such that $r < 1$, and put
\begin{equation}
	A(r) = \sum_{j=0}^\infty a_j \, r^j,
	\quad B(r) = \sum_{l=0}^\infty b_l \, r^l,
	\quad C(r) = \sum_{n=0}^\infty c_n \, r^n.
\end{equation}
If the $a_j$'s and $b_l$'s are bounded, for instance, then the $c_n$'s
grow at most linearly, and the series in the definitions of $A(r)$,
$B(r)$, and $C(r)$ converge absolutely when $0 \le r < 1$.
We also have that
\begin{equation}
	C(r) = A(r) \, B(r)
\end{equation}
for $0 \le r < 1$, because the series defining $C(r)$ is the Cauchy
product of the series defining $A(r)$ and $B(r)$.

	By definition, \emph{Abel summability} of $\sum_{j=0}^\infty
a_j$, $\sum_{l=0}^\infty b_l$, $\sum_{n=0}^\infty c_n$ means the
existence of the limits of $A(r)$, $B(r)$, $C(r)$ as $r \to 1$, $0 \le
r < 1$, in which case the Abel sum is defined to be the limit.
Ordinary convergence of an infinite series implies Abel summability,
with the Abel sum equal to the sum as the limit of the partial sums.
The Abel sum of the Cauchy product is equal to the product of the Abel
sums when they exist, and it follows that if the series converge, then
the sum of the Cauchy product is equal to the product of the sums of
the other two series.  See \cite{Rudin2} for more information.

	In the $p$-adic case the situation is simpler.
The series $\sum_{j=0}^\infty a_j$, $\sum_{l=0}^\infty b_l$
converge if and only if
\begin{equation}
	\lim_{j \to \infty} a_j = \lim_{l \to \infty} b_l = 0,
\end{equation}
and in this event
\begin{equation}
	\lim_{n \to \infty} c_n = 0,
\end{equation}
as one can check.  To see that the sum of the $c_n$'s is equal to the
product of the sums of the $a_j$'s and $b_l$'s, one can approximate by
finite sums and show that the relevant error terms are small.  This is
analogous to one of the arguments for absolutely convergent series of
real or complex numbers.

	Let us now consider the formal power series expansion
for the exponential function,
\begin{equation}
	E(z) = \sum_{n=0}^\infty \frac{z^n}{n!}.
\end{equation}
As usual, $n!$ denotes ``$n$ factorial'', which is the product of the
integers from $1$ to $n$, and which is interpreted as being equal to
$1$ when $n = 0$.  Formally we have that
\begin{equation}
\label{E(z + w) = E(z) E(w)}
	E(z + w) = E(z) \, E(w),
\end{equation}
in the sense that if one expands the series and collect terms
then the coefficients match up, as a result of the binomial theorem.

	Let us focus first on the case of complex numbers.  For each
$z \in {\bf C}$ one can show that the series defining $E(z)$ converges
absolutely, and indeed the terms tend to $0$ faster than a geometric
series.  In other words, the power series defining $E(z)$ has infinite
radius of convergence.  Hence the formal identity (\ref{E(z + w) =
E(z) E(w)}) does work for the actual sums for all $z, w \in {\bf C}$,
as a special case of the earlier discussion of Cauchy products.  In
particular, $E(z) \ne 0$ for all complex numbers $z$, with $1/E(z) =
E(-z)$.

	If $x$ is a nonnegative real number, then $E(x)$ is real and
$E(x) \ge 0$.  If $x$ is a real number and $x \le 0$, then $E(x)$ is a
real number such that $0 < E(x) \le 1$, since $E(x) = 1/E(-x)$.  If
$z$ is a complex number, then the complex conjugate of $E(z)$ is equal
to $E(\overline{z})$.  This follows from the series expansion for
$E(z)$, since the coefficients are real numbers.  As a consequence we
get that
\begin{equation}
	|E(z)|^2 = E(z) \, E(\overline{z})
		= E(z + \overline{z}) = E(2 \re z)
\end{equation}
for every complex number $z$.

	For the $p$-adic case we should begin by considering the
number of factors of $p$ in $n!$.  The number of positive integers
less than or equal to $n$ which are divisible by $p$ is equal to the
integer part of $n/p$.  For each positive integer $k$, the number of
positive integers less than or equal to $n$ which are divisible by
$p^k$ is equal to the integer part of $n / p_k$.  The total number of
factors of $p$ in $n!$ is equal to the sum of the integer parts of $n
/ p^k$ over all positive integers $k$.  This sum is less than $n / (p
- 1)$, by comparison with a geometric series.

	It follows that the series for $E(x)$ converges in ${\bf Q}_p$
when $x \in {\bf Q}_p$ has $p$-adic absolute value less than
$p^{-1/(p-1)}$.  This discussion follows the one in \cite{Fernando}
starting on p112.  As explained very nicely there, this condition may
seem a bit strange, since $1/(p-1)$ is an integer only when $p = 2$,
but in fact one may wish to consider $E(x)$ on extensions of ${\bf
Q}_p$ where the extension of the $p$-adic absolute value has nonzero
values other than integer powers of $p$.  A related point is that
there can be Galois actions on such an extension which then interact
with the exponential in a nice way, just as complex conjugation does
in the complex case.  At any rate, because of the ultrametric property
for the $p$-adic absolute value, a disk around $0$ is closed under
addition, and one again has the identity that the exponential of a sum
is equal to the product of the corresponding exponentials.

\section{Normed vector spaces}
\label{normed vector spaces}
\setcounter{equation}{0}

	Let $V$ be a vector space over the real or complex numbers.
By a \emph{seminorm} on $V$ we mean a nonnegative real-valued
function $N(v)$ on $V$ such that
\begin{equation}
\label{N(alpha v) = |alpha| N(v)}
	N(\alpha \, v) = |\alpha| \, N(v)
\end{equation}
for all real or complex numbers $\alpha$, as appropriate,
and all $v \in V$, and such that
\begin{equation}
\label{N(v + w) le N(v) + N(w)}
	N(v + w) \le N(v) + N(w)
\end{equation}
for all $v, w \in V$.

	Recall that a subset $E$ of $V$ is said to be
\emph{convex} if for every pair of vectors $v, w \in E$
and every real number $t$ with $0 \le t \le 1$ we have that
\begin{equation}
	t \, v + (1 - t) \, w \in E.
\end{equation}
Under the homogeneity condition (\ref{N(alpha v) = |alpha| N(v)}), one
can check that the triangle inequality (\ref{N(v + w) le N(v) + N(w)})
holds if and only if
\begin{equation}
	\{v \in V : N(v) \le 1\}
\end{equation}
is a convex subset of $V$.

	If $N(v)$ is a seminorm on $V$, and if $N(v) > 0$ for all $v
\in V$ with $v \ne 0$, then we say that $N(v)$ is a \emph{norm} on
$V$.  In this event we get a metric on $V$ given by $N(v - w)$.

	As a special case, suppose that $\langle v, w \rangle$ is an
\emph{inner product} on $V$, or more precisely a \emph{hermitian inner
product} in the complex case.  This means that $\langle v, w \rangle$
is a real or complex-valued function, according to whether $V$ is a
real or complex vector space, defined for $v, w \in V$, such that
$\langle v, w \rangle$ is a linear function of $v$ for each $w \in V$,
\begin{equation}
	\langle w, v \rangle = \langle v, w \rangle
\end{equation}
when $V$ is a real vector space and
\begin{equation}
	\langle w, v \rangle = \overline{\langle v, w \rangle}
\end{equation}
when $V$ is a complex vector space, and $\langle v, v \rangle$ is a
nonnegative real number for all $v \in V$ which is equal to $0$ if and
only if $v = 0$.  If we put
\begin{equation}
	\|v\| = \langle v, v \rangle^{1/2},
\end{equation}
then the Cauchy--Schwarz inequality states that
\begin{equation}
	|\langle v, w \rangle| \le \|v\| \, \|w\|
\end{equation}
for all $v, w \in V$.  This can be verified using the fact that
$\langle v + \alpha \, w, v + \alpha \, w \rangle$ is a nonnegative
real number for all scalars $\alpha$.  Using the Cauchy--Schwarz
inequality one can check that
\begin{equation}
	\|v + w\|^2 \le (\|v\| + \|w\|)^2
\end{equation}
for all $v, w \in V$, so that $\|v\|$ does in fact define a norm on $V$.

	If $V = {\bf R}^n$, then the standard inner product on $V$ is
given by
\begin{equation}
	\langle v, w \rangle = \sum_{j=1}^n v_j \, w_j.
\end{equation}
If $V = {\bf C}^n$, then the standard Hermitian inner product is
defined by
\begin{equation}
	\langle v, w \rangle = \sum_{j=1}^n v_j \, \overline{w_j}.
\end{equation}
The associated norm is the same as $\|v\|_2$ defined next.

	Let $p$ be a real number with $1 \le p < \infty$,
and put
\begin{equation}
	\|v\|_p = \biggl(\sum_{j=1}^n |v_j|^p \biggr)^{1/p}
\end{equation}
for $v = (v_1, \ldots, v_n)$ in ${\bf R}^n$ or ${\bf C}^n$.  We can
extend this to $p = \infty$ by setting
\begin{equation}
	\|v\|_\infty = \max(|v_1|, \ldots, |v_n|).
\end{equation}
For $1 \le p \le \infty$ we have that $\|v\|_p$ satisfies the
homogeneity condition (\ref{N(alpha v) = |alpha| N(v)}), and is equal
to $0$ if and only if $v = 0$.  When $p = 1, \infty$ one can check the
triangle inequality directly from the definitions, and when $p = 2$
this follows from the preceding discussion about inner product spaces.
In general when $1 < p < \infty$ one can check that the closed unit
ball associated to $\|v\|_p$ is a convex set, and hence that $\|v\|_p$
defines a norm, using the convexity of the function $t^p$ defined on
the nonnegative real numbers.

	For $1 \le p < \infty$ we have that
\begin{equation}
	\|v\|_\infty \le \|v\|_p
\end{equation}
for all $v$ in ${\bf R}^n$ or ${\bf C}^n$, by inspection.  Using this
one can verify more generally that
\begin{equation}
	\|v\|_q \le \|v\|_p
\end{equation}
for all $v$ in ${\bf R}^n$ or ${\bf C}^n$ when $1 \le p \le q \le
\infty$.

	If $r$ is a real number with $r \ge 1$, then
\begin{equation}
	\biggl(\frac{1}{n} \sum_{j=1}^n x_j \biggr)^r
		\le \frac{1}{n} \sum_{j=1}^n x_j^r
\end{equation}
for all nonnegative real numbers $x_1, \ldots, x_n$, by the convexity
of the function $t^r$ on the nonnegative real numbers.
As a result,
\begin{equation}
	\|v\|_p \le n^{1/p - 1/q} \, \|v\|_q
\end{equation}
when $v$ is an element of ${\bf R}^n$ or ${\bf C}^n$ and $1 \le p \le
q < \infty$.  This also works with $q = \infty$, $1/q = 0$, by
inspection.

	Suppose that $N(v)$ is a seminorm on ${\bf R}^n$ or ${\bf
C}^n$.  One can check that $N(v)$ is bounded by a constant times the
Euclidean norm $\|v\|_2$, or any other $\|v\|_p$ if one prefers, where
the constant can be estimated in terms of the values of $N$ at the
standard basis vectors.

	On any real or complex vector space $V$, if $N(v)$
is a seminorm on $V$, then
\begin{equation}
	N(v) \le N(w) + N(v - w)
\end{equation}
for all $v, w \in V$.  Similarly,
\begin{equation}
	N(w) \le N(v) + N(v - w)
\end{equation}
for all $v, w \in V$, and therefore
\begin{equation}
	|N(v) - N(w)| \le N(v - w)
\end{equation}
for all $v, w \in V$.  If $N$ is a norm, then $N$ is continuous with
respect to the metric associated to $N$.  If $V$ is ${\bf R}^n$ or
${\bf C}^n$, then $N$ is also continuous with respect to the standard
Euclidean metric, using the remark in the previous paragraph.  It
follows that there is a positive real number $\eta$ such that $N(v)
\ge \eta$ when $\|v\|_2 = 1$, which is to say that $v$ lies on the
standard Euclidean sphere, since $N(v)$ is a positive continuous
function on the sphere and the sphere is compact.

	One could define the notion of a norm just as well on a vector
space over a subfield of the real or complex numbers, like the
rational numbers.  One should be a bit careful, in that for instance
if $\alpha$ is an irrational number, then $N(x) = |x_1 - \alpha \,
x_2|$ defines a norm on ${\bf Q}^2$ which is more degenerate than
norms on ${\bf R}^2$ or ${\bf C}^2$.

	Instead one might consider vector spaces defined over the
rational or $p$-adic numbers with respect to the $p$-adic absolute
value function on scalars.  In this case one might consider the usual
triangle inequality for seminorms, as above, or the stronger
``ultrametric'' version requiring that the seminorm applied to a sum
of two vectors is less than or equal to the maximum of the values of
the seminorm at the two vectors.  Assuming homogeneity, this stronger
ultrametric version of the triangle inequality is equivalent to saying
that the set of vectors in the space with seminorm less than or equal
to $1$ is closed under addition.

	More generally one might consider vector spaces over fields
with absolute value functions, including extensions of the $p$-adic
numbers, as in \cite{Cassels1, Fernando}.  This may involve fields
which are not locally compact, and sometimes one is interested in
completeness instead, as in \cite{Cassels1, Fernando}.

\section{Dual spaces}
\label{dual spaces}
\setcounter{equation}{0}

	Let $V$ be a finite-dimensional real or complex vector space,
and let $V^*$ denote the dual vector space of linear functionals on
$V$.  Thus $V^*$ consists of the linear mappings from $V$ into the
real or complex numbers, whichever is the scalar field for $V$.  One
can add linear functionals and multiply them by scalars, so that $V^*$
is indeed a vector space with the same field of scalars as $V$.

	Suppose that $v_1, \ldots, v_n$ is a basis for $V$, so that
any vector $v \in V$ can be expressed in a unique way as a linear
combination of the $v_j$'s.  If $\lambda$ is a linear functional on
$V$, then $\lambda$ is uniquely determined by $\lambda(v_1), \ldots,
\lambda(v_n)$, since $\lambda(v)$ for any $v \in V$ can be computed
from the knowledge of these quantities and the coefficients of $v$ in
the basis.  Furthermore, for any collection of $n$ scalars $\alpha_1,
\ldots, \alpha_n$, there is a linear functional $\lambda$ on $V$ such
that $\lambda(v_j) = \alpha_j$ for each $j$.  In particular, the
dimension of $V^*$ is equal to the dimension of $V$.

	Now suppose that $V$ is also equipped with a norm $\|v\|$.
Let $\lambda$ be any linear functional on $V$, and put
\begin{equation}
	\|\lambda\|_* = \sup \{|\lambda(v)| : v \in V, \|v\| \le 1\}.
\end{equation}
To see that this is finite one can use an isomorphism between $V$ and
${\bf R}^n$ or ${\bf C}^n$, as appropriate, and the fact that $\|v\|$
is equivalent to a standard norm given explicitly in terms of
coordinates of vectors, as in Section \ref{normed vector spaces}.

	Equivalently, $\|\lambda\|_*$ can be characterized as a
nonnegative real number such that
\begin{equation}
	|\lambda(v)| \le \|\lambda\|_* \, \|v\|
\end{equation}
for all $v \in V$ and $\|\lambda\|_*$ is as small as possible.  One
can check that $\|\lambda\|_*$ defines a norm on the dual space $V^*$,
called the \emph{dual norm} associated to the norm $\|v\|$ on $V$.

	For instance, let $V$ be a real or complex vector space
equipped with an inner product $\langle v, w \rangle$.  For each $w
\in V$, $\lambda(v) = \langle v, w \rangle$ defines a linear
functional on $V$.  Using the Cauchy--Schwarz inequality one can check
that the dual norm of $\lambda$ is equal to the norm of $w$, with
respect to the norm on $V$ associated to the inner product.

	Now let $n$ be a positive integer, and let $V$ be ${\bf R}^n$
or ${\bf C}^n$.  A linear functional $\lambda$ on $V$ can be
represented explicitly as
\begin{equation}
\label{lambda(v) = sum_{j=1}^n v_j w_j}
	\lambda(v) = \sum_{j=1}^n v_j \, w_j,
\end{equation}
$v = (v_1, \ldots, v_n)$, where $w = (w_1, \ldots, w_n)$ is an element
of ${\bf R}^n$ or ${\bf C}^n$, as appropriate.  If $V$ is equipped
with the norm $\|v\|_1$ as in Section \ref{normed vector spaces}, then
we have that $|\lambda(v)| \le \|w\|_\infty \, \|v\|_1$ for all $v \in
V$, just by the triangle inequality.  One can also verify that
$\|w\|_\infty$ is the smallest nonnegative real number with this
property.  In other words, if we use the norm $\|\cdot \|_1$ on $V$,
then the dual norm is given by $\|\cdot \|_\infty$.

	Now suppose that we use the norm $\|v\|_\infty$ on $V$.  As in
the previous paragraph we have that $|\lambda(v)| \le \|w\|_1 \,
\|v\|_\infty$ for all $v \in V$ when $\lambda$ is associated to $w$ as
in (\ref{lambda(v) = sum_{j=1}^n v_j w_j}), by the triangle
inequality, and that $\|w\|_1$ is the smallest nonnegative real number
with this property, so that the dual norm of $\lambda$ is exactly
$\|w\|_1$.

	Let $p$, $q$ be real numbers with $1 < p, q < \infty$ and $1/p
+ 1/q = 1$.  In this case we say that $p$, $q$ are \emph{conjugate
exponents}.  One can check that
\begin{equation}
	a \, b \le \frac{a^p}{p} + \frac{b^q}{q}
\end{equation}
for any nonnegative real numbers $a$, $b$, and indeed one can view
this as a consequence of the convexity of the exponential function on
the real line.  If $a_1, \ldots, a_n$ and $b_1, \ldots, b_n$ are
nonnegative real numbers, then \emph{H\"older's inequality} states
that
\begin{equation}
	\sum_{j=1}^n a_j \, b_j \le \biggl(\sum_{k=1}^n a_k^p \biggr)^{1/p}
		\, \biggl(\sum_{l=1}^n b_l^q \biggr)^{1/q}.
\end{equation}
This follows from the previous inequality when $\sum_k a_k^p \le 1$
and $\sum_l b_l^q \le 1$, just by applying the inequality termwise and
summing, and one can derive the general case from this by a scaling
argument.

	If we use the norm $\|v\|_p$ on $V$, and if the linear
function $\lambda$ is associated to an $n$-tuple $w$ as before, then
we have that $|\lambda(v)| \le \|w\|_q \, \|v\|_p$ for all $v \in V$
by H\"older's inequality.  For a fixed $w$ one can choose $v \ne 0$ so
that this inequality becomes an equality, as one can check.  As a
result, the dual norm of $\lambda$ associated to the norm $\|v\|_p$ on
$V$ is equal to $\|w\|_q$.

\section{Operator norms}
\label{operator norms}
\setcounter{equation}{0}

	Let $V$ be a finite-dimensional real or complex vector space,
and let $\mathcal{L}(V)$ denote the collection of linear mappings from
$V$ into itself.  Thus $\mathcal{L}(V)$ is a vector space in a natural
way, since one can add linear transformations on $V$ and one can
multiply them by scalars.  Moreover one can compose linear
transformations on $V$, which provides a kind of multiplication on
$\mathcal{L}(V)$, making it an algebra rather than simply a vector
space.  The identity transformation $I$ on $V$, which sends each
vector $v \in V$ to itself, acts as a multiplicative identity element
in the algebra, since the composition of any linear transformation $T$
on $V$ with $I$ is equal to $T$.

	Let $v_1, \ldots, v_n$ be a basis for $V$.  If $T$ is a linear
transformation on $V$, then $T$ is uniquely determined by its values
on the $v_j$'s.  The image of each $v_j$ under $T$ is a vector in $V$
and therefore characterized by its $n$ coefficients with respect to
the basis $v_1, \ldots, v_n$.  Conversely one can start with $n^2$
scalars, which can be arranged naturally into an $n \times n$ matrix,
and get a linear transformation $T$ on $V$ for which the given scalars
are the coefficients of the $T(v_j)$'s in the basis.  In particular
$\mathcal{L}(V)$ has dimension $n^2$ as a vector space.

	Suppose that $V$ is equipped with a norm $\|v\|$.
If $T$ is a linear transformation on $V$, then put
\begin{equation}
	\|T\|_{op} = \sup \{\|T(v)\| : v \in V, \|v\| \le 1\}.
\end{equation}
That this is finite can be seen using an isomorphism between $V$ and
${\bf R}^n$ or ${\bf C}^n$, as appropriate, and the equivalence of
$\|v\|$ with a standard norm given in terms of coordinates.

	One can also characterize $\|T\|_{op}$, called the
\emph{operator norm} of $T$ associated to the norm $\|v\|$
on $V$, as the smallest nonnegative real number such that
\begin{equation}
	\|T(v)\| \le \|T\|_* \, \|v\|
\end{equation}
for all $v \in V$.  It is easy to verify that this does indeed define
a norm on the vector space of linear transformations on $V$, and it
enjoys the additional property that
\begin{equation}
	\|T_1 \circ T_2\|_{op}
		\le \|T_1\|_{op} \, \|T_2\|_{op}
\end{equation}
for any linear operators $T_1$, $T_2$ on $V$.

	Of course the norm of the identity transformation $I$ on $V$
is equal to $1$.  If $T$ is an invertible linear transformation on
$V$, so that there is a linear transformation $T^{-1}$ on $V$ whose
composition with $T$ is equal to $I$, then
\begin{equation}
	1 = \|I\|_{op} \le \|T\|_{op} \, \|T^{-1}\|_{op}.
\end{equation}
Suppose that $A$ is a linear transformation on $V$, $v$ is a nonzero
vector in $V$, and that $\alpha$ is a real or complex number, as
appropriate.  We say that $v$ is an \emph{eigenvector} for $A$ with
\emph{eigenvalue} $\alpha$ if
\begin{equation}
	A(v) = \alpha \, v.
\end{equation}
In this event
\begin{equation}
	|\alpha| \le \|A\|_{op}.
\end{equation}

	Let $V$ be ${\bf R}^n$ or ${\bf C}^n$, and let $e_1, \ldots,
e_n$ denote the standard basis vectors for $V$, which is to say that
the $l$th component of $e_j$ is equal to $1$ when $l = j$ and is equal
to $0$ otherwise.  Suppose that we use the norm $\|v\|_1$ from Section
\ref{normed vector spaces} for $V$.  If $T$ is a linear operator on
$V$, then the operator norm of $T$ with respect to this norm on $V$ is
equal to the maximum of the norms of $T(e_1), \ldots, T(e_n)$.  This
is not difficult to verify just from the definitions.  Of course this
can be expressed explicitly in terms of the absolute values of the
entries of the matrix associated to $T$ with respect to the standard
basis of $e_j$'s.

	Suppose instead that we use the norm $\|v\|_\infty$ from
Section \ref{normed vector spaces}.  We can think of $T$ as being
described by $n$ linear functionals $\lambda_1, \ldots, \lambda_n$ on
$V$, where $\lambda_j(v)$ is the same as the $j$th component of $T(v)$
for all $v \in V$.  The operator norm of $T$ with respect to the norm
$\|v\|_\infty$ on $V$ is equal to the maximum of the dual norms of
$\lambda_1, \ldots, \lambda_n$, as one can easily verify.  The dual
norm associated to $\|v\|_\infty$ was determined in the previous
section, and thus the operator norm of $T$ can again be given
explicitly in terms of the absolute values of the matrix entries of
$T$ with respect to the standard basis in this case.

	Let $T$ be a linear operator on $V$ whose operator norm with
respect to each of $\|\cdot \|_1$ and $\|\cdot \|_\infty$ is less than
or equal to $1$.  In other words, assume that
\begin{equation}
	\|T(v)\|_1 \le \|v\|_1
\end{equation}
and that
\begin{equation}
	\|T(v)\|_\infty \le \|v\|_\infty
\end{equation}
for all $v \in V$.  In terms of the matrix of $T$ associated to the
standard basis $e_1, \ldots, e_n$, this is equivalent to saying that
the sum of the absolute values of the matrix entries in any row or
column is less than or equal to $1$.  A result of Schur implies that
\begin{equation}
	\|T(v)\|_p \le \|v\|_p
\end{equation}
for any $p$, $1 < p < \infty$, and all $v \in V$, which is to say that
the operator norm of $T$ with respect to $\|\cdot \|_p$ is also less
than or equal to $1$.

	To show this we may as well assume that the matrix entries of
$T$ are nonnegative real numbers, and we may as well restrict our
attention to vectors $v$ whose components are nonnegative real
numbers.  In other words, we can reduce to this case by putting in
absolute values everywhere and applying the triangle inequality
repeatedly.  The hypotheses on $T$ still hold if we replace $T$ with
the linear transformation whose matrix entries are the absolute values
of the matrix entries of $T$.

	Under these conditions, one can check that the $p$th power of
the $j$th component of $T(v)$ is less than or equal to the $j$th
component of $T$ applied to the vector given by the $p$th power of the
components of $v$.  This follows from the convexity of the function
$t^p$ on the nonnegative real numbers, using the fact that the
operator norm of $T$ with respect to $\|\cdot \|_\infty$ is less than
or equal to $1$.  Because $T$ has operator norm less than or equal to
$1$ with respect to $\|\cdot \|_1$, it follows that the sum of the
$p$th powers of the components of $T(v)$ is less than or equal to the
sum of the $p$th powers of the components of $v$.  This says exactly
that $\|T(v)\|_p^p \le \|v\|_p^p$, as desired.

	As another special case, let $T$ be a linear operator on ${\bf
R}^n$ or ${\bf C}^n$ which is diagonalized by the standard basis $e_1,
\ldots, e_n$.  That is, we assume that there are real or complex
numbers $\alpha_1, \ldots, \alpha_n$, as appropriate,
so that
\begin{equation}
	T(e_j) = \alpha_j \, e_j
\end{equation}
for $j = 1, \ldots, n$.  In this event the operator norm of $T$ is
equal to
\begin{equation}
	\max (|\alpha_1|, \ldots, |\alpha_n|)
\end{equation}
with respect to any of the norms $\|v\|_p$, $1 \le p \le \infty$.

	Now suppose that $V$ is a finite-dimensional real or complex
vector space equipped with an inner product $\langle v, w \rangle$,
and let $A$ be a linear transformation on $V$.  We say that $A$ is
\emph{self-adjoint} if
\begin{equation}
	\langle A(v), w \rangle = \langle v, A(w) \rangle
\end{equation}
for all $v, w \in V$.  A famous theorem states that $A$ can be
diagonalized in an orthonormal basis in this situation.  In other
words, there exist vectors $v_1, \ldots, v_n in V$ and real numbers
$\alpha_1, \ldots, \alpha_n$ such that
\begin{equation}
	\langle v_j, v_l \rangle = 0
\end{equation}
when $j \ne l$, $\langle v_j, v_j \rangle = 1$ for all $j$, every
element of $V$ can be expressed as a linear combination of the
$v_j$'s, and $A(v_j) = \alpha_j \, v_j$ for each $j$.  The operator
norm of $A$ is then equal to the maximum of $|\alpha_1|, \ldots,
|\alpha_n|$ with respect to the norm associated to the inner product.

\section{Geometry of numbers}
\label{geometry of numbers}
\setcounter{equation}{0}

	Consider ${\bf Q}^n$, the subset of ${\bf R}^n$ consisting of
points with rational coordinates, and suppose that that $N(v)$ is a
seminorm on ${\bf Q}^n$.  As in Section \ref{normed vector spaces},
this means that $N(v)$ is a nonnegative real-valued function defined
for $v \in {\bf Q}^n$ such that $N(\alpha \, v) = |\alpha| \, N(v)$
for all $\alpha \in {\bf Q}$ and $v \in V$, and $N(v + w) \le N(v) +
N(w)$ for all $v, w \in V$.  Because $N(v) \le N(w) + N(v-w)$ and
$N(w) \le N(v) + N(v-w)$ for all $v, w \in {\bf Q}^n$, we have that
\begin{equation}
	|N(v) - N(w)| \le N(v-w)
\end{equation}
for all $v, w \in {\bf Q}^n$, as before.  We also have that $N(v)$ is
bounded by a constant multiple of the Euclidean norm on ${\bf Q}^n$,
with an estimate in terms of $N(e_j)$, $1 \le j \le n$, where the
$e_j$'s are the standard basis vectors in ${\bf R}^n$.  It follows
that $N(v)$ is a uniformly continuous function on ${\bf Q}^n$, and
therefore has a unique continuous extension to a function on ${\bf
R}^n$ which is in fact a seminorm on ${\bf R}^n$.

	The extension of $N$ to ${\bf R}^n$ may or may not be a norm,
even if $N$ is a norm on ${\bf Q}^n$.  For if $\theta$ is any real
number, then $N(v) = |v_1 - \theta \, v_2|$ defines a seminorm on
${\bf R}^2$ and on ${\bf Q}^2$ by restriction.  As in Section
\ref{normed vector spaces}, if $\theta$ is irrational, then $N(v)$ is
a norm on ${\bf Q}^2$, but it is not a norm on ${\bf R}^2$ for any
$\theta$.  At any rate, if we start with a norm on ${\bf Q}^n$ for
some $n$, we can extend it to a seminorm on ${\bf R}^n$, and it is
interesting to consider the interplay between the norm on ${\bf Q}^n$
/ seminorm on ${\bf R}^n$ and arithmetic.

	Suppose now that $U$ is an open subset of ${\bf R}^n$, and let
${\bf Z}^n$ denote the subset of ${\bf R}^n$ of points with integer
coordinates, which is of course closed under addition.  If the volume
of $U$ is strictly larger than $1$, then there are points $x, y \in U$
with $x \ne y$ and $x - y \in {\bf Z}^n$.  To see this it is
convenient to think of the quotient of ${\bf R}^n$ by ${\bf Z}^n$,
which is a torus whose total volume is equal to $1$, and the natural
projection $p$ from ${\bf R}^n$ onto the quotient.  The existence of
distinct points $x$, $y$ in $U$ whose difference is an element of
${\bf Z}^n$ is equivalent to saying that the restriction of $p$ is not
one-to-one, which follows immediately if the volume of $U$ is strictly
larger than the volume of the quotient torus, which is equal to $1$.

	Assume further that $U$ is symmetric about the origin, so that
$w \in U$ implies $-w \in U$.  We can rephrase the previous conclusion
then to say that there are points $x, y \in U$ such that $x + y$ is a
nonzero element of ${\bf Z}^n$.  If $U$ is also convex, then $(x +
y)/2$ is a nonzero element of $U$ which lies in $(1/2) {\bf Z}^n$,
which is to say that its coordinates are integers or half-integers.
We can rephrase this again by saying that if $U$ is a convex open
subset of ${\bf R}^n$ which is symmetric about the origin and which
has volume strictly larger than $2^n$, then $U$ contains a nonzero
element of ${\bf Z}^n$.  This is a version of the basic existence
result in the geometry of numbers.

\section{Linear groups}
\label{linear groups}
\setcounter{equation}{0}

	Let $V$ be a finite-dimensional real or complex vector space,
and let $GL(V)$ denote the group of invertible linear transformations
on $V$.  Of course any nonzero multiple of the identity operator $I$
is invertible.  As a subset of the vector space $\mathcal{L}(V)$ of
linear transformations on $V$, $GL(V)$ is open, since it consists
simply of the linear transformations with nonzero determinant.

	We can also look at this in terms of norms.  Let $\|v\|$ be a
norm on $V$, so that we have an associated norm $\|T\|_{op}$ for
linear operators on $V$.  If $T$ is an invertible linear
transformation on $V$, then there is a real number $c > 0$ such that
\begin{equation}
	c \, \|v\| \le \|T(v)\|
\end{equation}
for all $v \in V$, namely, $c = 1/\|T^{-1}\|_{op}$.  If $A$ is a
linear operator on $V$ such that $\|A\|_{op} < c$, then
\begin{equation}
	(c - \|A\|_{op}) \, \|v\| \le \|(T + A)(v)\|
\end{equation}
for all $v \in V$.  It follows that $T + A$ is injective, and hence
invertible, since $V$ is assumed to have finite dimension.

	Once we specify a norm $\|v\|$ on $V$, we get a nice subgroup
of $GL(V)$, namely the group of linear isometries on $V$, which are
the linear mappings $T$ from $V$ to itself such that
\begin{equation}
	\|T(v)\| = \|v\|
\end{equation}
for all $v \in V$.  This is the same as saying that both $T$ and
$T^{-1}$ have norm equal to $1$, and of course the identity operator
$I$ is always an isometry.  The group of isometries on $V$ is a
compact subset of $GL(V)$, because it is closed and bounded.  If the
norm on $V$ comes from an inner product, then the group of isometries
is quite rich, and is known as an orthogonal or unitary group,
according to whether $V$ is a real or complex vector space.

	Let us take $V = {\bf R}^n$, and consider the group of
invertible linear transformations $T$ on ${\bf R}^n$ which map ${\bf
Z}^n$ onto itself.  A linear mapping $T$ on ${\bf R}^n$ maps ${\bf
Z}^n$ into itself if and only if the matrix associated to $T$ and the
standard basis $e_1, \ldots, e_n$ in ${\bf R}^n$ has integer entries.
In order that $T$ be an invertible linear transformation on ${\bf
R}^n$ which takes ${\bf Z}^n$ onto itself the matrices associated to
both $T$ and $T^{-1}$ should have integer entries.  This is equivalent
to saying that the matrix associated to $T$ has integer entries and
determinant equal to $\pm 1$.

	Now consider ${\bf Q}_p^n$, the space of $n$-tuples of
$p$-adic numbers, as a vector space over ${\bf Q}_p$ with respect to
coordinatewise addition and scalar multiplication.  The group of
invertible linear transformations on ${\bf Q}_p^n$ is described by the
condition that the determinant is nonzero, and is an open subset of
the vector space of all linear transformations on ${\bf Q}_p$, which
can be identified with ${\bf Q}_p^{n^2}$.  A natural norm on ${\bf
Q}_p^n$ is given by
\begin{equation}
	N(v) = \max (|v_1|_p, \ldots, |v_n|_p)
\end{equation}
for $v = (v_1, \ldots, v_n) \in {\bf Q}_p^n$.  The linear mappings $T$
on $V$ which are isometries with respect to this norm can be
characterized by the condition that the matrix with respect to the
standard basis has entries in ${\bf Z}_p$ and the determinant has
$p$-adic absolute value equal to $1$, so that the inverse matrix also
has entries in ${\bf Z}_p$.  This is a compact subgroup of the group
of all invertible linear transformations on ${\bf Q}_p$.

\section{Trace norms}
\label{trace norms}
\setcounter{equation}{0}

	Let $V$ be a finite-dimensional real or complex vector space
equipped with an inner product $\langle v, w \rangle$.  As before, a
linear transformation $A$ on $V$ is self-adjoint if
\begin{equation}
	\langle A(v), w \rangle = \langle v, A(w) \rangle
\end{equation}
for all $v, w \in V$.  In this case $A$ can be diagonalized in an
orthonormal basis for $V$, which is to say that there is an
orthonormal basis $v_1, \ldots, v_n$ for $V$ and real numbers
$\alpha_1, \ldots, \alpha_n$ such that $A(v_j) = \alpha_j \, v_j$ for
$j = 1, \ldots, n$.  The $\alpha_j$'s are the eigenvalues of $A$, and
it is easy to verify directly that the eigenvalues of $A$ are real
numbers even if $V$ is a complex vector space.  Notice that the
self-adjoint linear operators on $V$ form a real vector space, which
is to say that the sum of two self-adjoint linear operators on $V$ is
a self-adjoint linear operator on $V$ and that the product of a real
number and a self-adjoint linear operator on $V$ is a self-adjoint
linear operator on $V$, and that one should use real scalars for this
even if $V$ is a complex vector space.

	For $A$ as in the previous paragraph and $p \in {\bf R}$ such
that $1 \le p < \infty$, put
\begin{equation}
	\|A\|_{\mathcal{S}_p} = \biggl(\sum_{j=1}^n |\alpha_j|^p \biggr)^{1/p}.
\end{equation}
This is the \emph{Schatten $p$-class} or $\mathcal{S}_p$ norm of $A$,
although the fact that the triangle inequality holds for the
$\mathcal{S}_p$ norm is somewhat tricky and will be discussed further
in a moment.  When $p = \infty$ one can define the
$\mathcal{S}_\infty$ norm of $A$ to be the maximum of the
$|\alpha_j|'s$, which is the same as the ordinary operator norm of
$A$.

	Let $w_1, \ldots, w_n$ be another orthonormal basis for $V$,
in addition to the orthonormal basis $v_1, \ldots, v_n$ of
eigenvectors for $A$, and consider
\begin{equation}
\label{(sum_{l=1}^n |langle A(w_l), w_l rangle|^p )^{1/p}}
	\biggl(\sum_{l=1}^n |\langle A(w_l), w_l \rangle|^p \biggr)^{1/p}.
\end{equation}
This is equal to $\|A\|_{\mathcal{S}_p}$ when $w_j = v_j$ for each
$j$.  In general,
\begin{equation}
	\langle A(w_l), w_l \rangle
		\sum_{j=1}^n \alpha_j \, |\langle v_j, w_l \rangle|^2,
\end{equation}
since $A(w) = \sum_{j=1}^n \alpha_j \, \langle w, v_j \rangle \, v_j$
for all $w \in V$.

	Because $v_1, \ldots, v_n$ is an orthonormal basis,
$\sum_{j=1}^n |\langle v_j, w \rangle|^2 = \|w\|^2$ for all $w \in V$.
Similarly, $\sum_{l=1}^n |\langle v, w_l \rangle|^2 = \|v\|^2$ for all
$v \in V$.  It follows that the sum of $|\langle v_j, w_l \rangle|^2$
over $j$ for a fixed $l$, or the sum over $l$ for a fixed $j$, is
equal to $1$.

	By the result of Schur mentioned in Section \ref{operator
norms} we get that
\begin{equation}
  \biggl(\sum_{l=1}^n 
\biggl|\sum_{j=1}^n x_j \, |\langle v_j, w_l \rangle|^2 \biggr|^p \biggr)^{1/p}
	\le \biggl(\sum_{j=1}^n |x_j|^p \biggr)^{1/p}
\end{equation}
for all real or complex numbers $x_1, \ldots, x_n$ and $1 \le p <
\infty$.  This implies that (\ref{(sum_{l=1}^n |langle A(w_l), w_l
rangle|^p )^{1/p}}) is always less than or equal to
$\|A\|_{\mathcal{S}_p}$ for any orthonormal basis $w_1, \ldots, w_n$
on $V$.  Therefore $\|A\|_{\mathcal{S}_p}$ is the same as the maximum
of (\ref{(sum_{l=1}^n |langle A(w_l), w_l rangle|^p )^{1/p}}) over all
orthonormal bases $w_1, \ldots, w_n$ on $V$, and it follows that
$\|A\|_{\mathcal{S}_p}$ is indeed a norm on the real vector space of
self-adjoint linear operators on $V$.

\section{Vector spaces}
\label{vector spaces}
\setcounter{equation}{0}

	Let $k$ be a field, and let $V$ be a vector space over $k$.
If $v_1, \ldots, v_n$ are elements of $V$, then their \emph{span} is
denoted $\span (v_1, \ldots, v_n)$ and consists of all linear
combinations
\begin{equation}
\label{alpha_1 v_1 + cdots + alpha_n v_n}
	\alpha_1 \, v_1 + \cdots + \alpha_n \, v_n
\end{equation}
of $v_1, \ldots, v_n$, with $\alpha_1, \ldots, \alpha_n \in k$.  The
span of $v_1, \ldots, v_n$ is a linear subspace of $V$, which is to
say that it is closed under addition and scalar multiplication.

	A collection of vectors $v_1, \ldots, v_n$ in $V$ is said to
be \emph{linearly independent} if the linear combination (\ref{alpha_1
v_1 + cdots + alpha_n v_n}) of $v_1, \ldots, v_n$ is equal to $0$ if
and only if the scalars $\alpha_j$ are all equal to $0$, $1 \le j \le
n$.  This is equivalent to saying that each element of the span of
$v_1, \ldots, v_n$ can be expressed as a linear combination of $v_1,
\ldots, v_n$ in a unique manner.  A collection of vectors $v_1,
\ldots, v_n$ in $V$ is said to be \emph{linearly dependent} if the
$v_j$'s are not linearly independent.  This is equivalent to saying
that one of the $v_j$'s lies in the span of the others.

	Suppose that $v_1, \ldots, v_n$ and $w_1, \ldots, w_m$ are
vectors in $V$, with each $v_j$ an element of the span of $w_1,
\ldots, w_m$.  If $n > m$, then the $v_j$'s are linearly dependent.
Basically this reduces to the fact that a homogeneous system of $m$
linear equations and $n$ variables has a nontrivial solution when $n >
m$.

	We say that the vector space $V$ has finite dimension if there
is a finite collection of vectors in $V$ which span $V$.  The smallest
number of vectors in $V$ needed to span $V$ is called the dimension of
$V$.  If $V$ is the span of $v_1, \ldots, v_n$ and $n$ is as small as
possible, then $v_1, \ldots, v_n$ are linearly independent.  A
collection of vectors $v_1, \ldots, v_n$ in $V$ which are linearly
independent and whose span is equal to $V$ is called a basis for $V$.

\section{Algebras}
\label{algebras}
\setcounter{equation}{0}

	Let $k$ be a field, and let $\mathcal{A}$ be an algebra over
$k$.  This means that $\mathcal{A}$ is a vector space over $k$
equipped with a binary operation of multiplication which is
associative and which is bilinear with respect to the vector space
operations on $\mathcal{A}$.  We shall also assume that $\mathcal{A}$
has a nonzero multiplicative identity element $e$, and hence that
$\mathcal{A}$ contains a copy of $k$.

	Of course $k$ itself is a $1$-dimensional algebra over $k$.
If $E$ is a nonempty set, then the vector space $\mathcal{F}(E, k)$ of
$k$-valued functions on $E$ is an algebra over $k$, with the constant
function equal to $1$ at each point as the identity element.  If $V$
is a vector space over $k$, $V \ne \{0\}$, then the algebra
$\mathcal{L}(V)$ of linear operators on $V$ is an algebra over $k$,
with the identity operator $I$ as the multiplicative identity element.
Notice that for any algebra $\mathcal{A}$ over $k$ with nonzero
identity element, one can embed $\mathcal{A}$ into the algebra of
linear operators on $V = \mathcal{A}$, viewed simply as a vector
space.  Namely, each element $x$ of $\mathcal{A}$ induces a linear
transformation on $\mathcal{A}$ given by left multiplication by $x$.

	An algebra $\mathcal{A}$ over $k$ is said to be
\emph{commutative} if the operation of multiplication on $\mathcal{A}$
is commutative, i.e., if $x \, y = y \, x$ for all $x, y \in
\mathcal{A}$.  Let $\mathcal{A}$ be an algebra over $k$ with nonzero
multiplicative identity element $e$, which may or may not be
commutative, and let $x$ be any element of $\mathcal{A}$.  One can get
a subalgebra of $\mathcal{A}$ consisting of linear combinations of
powers of $x$ and $e$, and this subalgebra is automatically
commutative.

\section{Eigenvalues}
\label{eigenvalues}
\setcounter{equation}{0}

	Let $k$ be a field and let $V$ be a vector space over $k$ with
positive finite dimension equal to $n$.  Suppose that $A$ is a linear
transformation on $V$.  An element $\alpha$ of $k$ is said to be an
\emph{eigenvalue} for $A$ if there is a nonzero \emph{eigenvector}
corresponding to $A$, i.e., a vector $v \in V$, $v \ne 0$, such that
\begin{equation}
	A(v) = \alpha \, v.
\end{equation}
In this case $A - \alpha \, I$ has nontrivial kernel, where $I$ is the
identity operator on $V$, and hence is not invertible.  Conversely, if
$A - \alpha \, I$ is not invertible for some $\alpha \in k$, then $A -
\alpha \, I$ has nontrivial kernel, since $V$ is assumed to be
finite-dimensional, and this implies that $\alpha$ is an eigenvalue
for $A$.

	As a vector space over $k$, the algebra $\mathcal{L}(V)$ of
linear transformations on $V$ has dimension equal to $n^2$.  If $T$ is
a linear transformation on $V$, then there is a positive integer $l
\le n^2$ such that $T^l$ can be expressed as a linear combination of
the $T^j$'s with $1 \le j < l$ and the identity operator $I$, since
otherwise $\mathcal{L}(V)$ would have dimension larger than $n^2$.
Here $T^j$ refers to the product of $j$ $T$'s when $j$ is a positive
integer.  The Cayley--Hamilton Theorem gives a more precise version of
this, with $T^n$ expressed as a linear combination of $T^j$, $1 \le j
< n$, and the identity operator.

	At any rate there is a positive integer $l$ and scalars $c_0,
\ldots, c_{l-1} \in k$ such that
\begin{equation}
	T^l = c_{l-1} \, T^{l-1} + \cdots + c_1 \, T + c_0 \, I.
\end{equation}
If $T$ is invertible and $l$ is as small as possible, then $c_0 \ne
0$, because otherwise we could remove a factor of $T$ and express
$T^{l-1}$ as a linear combination of smaller powers of $T$ and the
identity operator.  As a result, if $T$ is an invertible linear
transformation on $V$, then the inverse of $T$ can be expressed as a
linear combination of powers of $T$ and the identity operator on $V$.

	Let $A$ be a linear operator on $V$ again, and let
$\mathcal{A}$ denote the commutative subalgebra of $\mathcal{L}(V)$
consisting of linear operators on $V$ which can be expressed as a
linear combination of powers of $A$ and the identity transformation.
From the previous remarks it follows that if $T \in \mathcal{A}$ and
$T$ is invertible as a linear operator on $V$, then the inverse of $T$
also lies in $\mathcal{A}$.

	Suppose that $\alpha \in k$ is an eigenvalue of $A$, and let
$v$ be a nonzero vector in $V$ which is an eigenvector for $A$ with
eigenvalue $\alpha$.  Of course $v$ is also an eigenvector for the
identity operator with eigenvalue $1$, and $v$ is an eigenvector for
$A^j$ with eigenvalue $\alpha^j$ for every positive integer $j$.  If
$T$ is any element of $\mathcal{A}$, then $v$ is an eigenvector for
$T$ with eigenvalue that we shall denote $\phi(T)$.  One can check
that $\phi$ defines an algebra homomorphism from $\mathcal{A}$ onto
$k$.

	Conversely, suppose that we start with a homomorphism $\phi$
from $\mathcal{A}$ into $k$ which is not identically equal to $0$.  It
follows that $\phi(I) = 1$, and that $\phi$ maps $\mathcal{A}$ onto
$k$.  If $T$ is any element of $\mathcal{A}$ which is invertible as a
linear transformation on $V$, and hence has its inverse in
$\mathcal{A}$, then $\phi(T) \ne 0$.

	Put $\alpha = \phi(A)$.  It follows that $\phi$ applied to $A
- \alpha \, I$ is equal to $0$, and hence $A - \alpha \, I$ is not
invertible, which is to say that $\alpha$ is an eigenvalue of $A$.  In
short, the set of eigenvalues of $A$ is equal to the set of values of
nonzero homomorphisms from $\mathcal{A}$ into $k$ at $A$.  Notice also
that each nonzero homomorphism from $\mathcal{A}$ into $k$ is
determined by its value at $A$.  Therefore the number of nonzero
homomorphisms from $\mathcal{A}$ into $k$ is equal to the number of
distinct eigenvalues of $A$.

\section{Polynomials, formal power series}
\label{polynomials, formal power series}
\setcounter{equation}{0}

	Let $k$ be a field, and let us write $\Sigma_0(k)$,
$\Sigma(k)$ for the spaces of sequences $\{a_j\}_{j=0}^\infty$ with
$a_j \in k$ for all $j$ and with $a_j = 0$ for sufficiently large $j$,
depending on the sequence, in the case of $\Sigma(k)$.  Thus
$\Sigma_0(k)$, $\Sigma(k)$ are vector spaces over $k$ with respect to
termwise addition and scalar multiplcation, and $\Sigma_0(k)$ is a
linear subspace of $\Sigma(k)$.  If $\{a_j\}_{j=0}^\infty$,
$\{b_l\}_{l=0}^\infty$ are sequences in $\Sigma(k)$, then their Cauchy
product is the sequence $\{c_n\}_{n=0}^\infty$ in $\Sigma(k)$ defined
by
\begin{equation}
	c_n = \sum_{j=0}^n a_j \, b_{n-j}
\end{equation}
for all $n \ge 0$, and this lies in $\Sigma_0(k)$ if
$\{a_j\}_{j=0}^\infty$ and $\{b_l\}_{l=0}^\infty$ do.  With respect to
this product, $\Sigma_0(k)$ and $\Sigma(k)$ become commutative
algebras over $k$.

	Let us write $\mathcal{P}(k)$, $\mathcal{PS}(k)$ for the
algebras of polynomials and formal power series over $k$.
Thus an element of $\mathcal{P}(k)$ can be expressed as
\begin{equation}
\label{a_n t^n + a_{n-1} t^{n-1} + cdots + a_0}
	a_n \, t^n + a_{n-1} \, t^{n-1} + \cdots + a_0
\end{equation}
for some $c_0, \ldots, c_n \in k$, while an element of
$\mathcal{PS}(k)$ can be expressed as
\begin{equation}
	\sum_{j=0}^\infty a_j \, t^j
\end{equation}
for some sequence of coefficients $a_0, a_1, \ldots$ in $k$.  Here $t$
is an indeterminant, and one can add and multiply polynomials or power
series in the usual manner, so that $\mathcal{P}(k)$ and
$\mathcal{PS}(k)$ are commutative algebras over $k$.  Of course
$\mathcal{P}(k)$ is the subalgebra of $\mathcal{PS}(k)$ with all but
finitely many coefficients equal to $0$.

	There is a natural one-to-one correspondence between
$\Sigma(k)$ and $\mathcal{PS}(k)$, in which a sequence with terms in
$k$ is associated to the formal power series with that sequence of
coefficients.  In this correspondence $\Sigma_0(k)$ is mapped onto
$\mathcal{P}(k)$.  More precisely this defines an isomorphism between
$\Sigma_0(k)$ and $\mathcal{P}(k)$, and between $\Sigma(k)$ and
$\mathcal{PS}(k)$, as algebras over $k$, which is to say that addition
and multiplication are preserved.

	Suppose that $\mathcal{A}$ is any algebra over $k$ with
nonzero multiplicative identity element $e$, and let $x$ be an element
of $\mathcal{A}$.  If $p(t)$ is a polynomial over $k$, with $p(t)$
given by (\ref{a_n t^n + a_{n-1} t^{n-1} + cdots + a_0}) for some
$a_0, \ldots, a_n \in k$, then we can define $p(x)$ to be the element
of $\mathcal{A}$ given by
\begin{equation}
	a_n \, x^n + a_{n-1} \, x^{n-1} + \cdots + a_0 \, e.
\end{equation}
Notice that if $p_1(t)$, $p_2(t)$ are polynomials over $k$, then $(p_1
+ p_2)(x) = p_1(x) + p_2(x)$ and $(p_1 \, p_2)(x) = p_1(x) \, p_2(x)$.
As a special case, we can take $\mathcal{A}$ equal to $k$, viewed as a
one-dimensional algebra over itself, and then we are simply saying
that a formal polynomial $p(t) \in \mathcal{P}(k)$ defines a function
on $k$ in the usual manner.  As another special case, we can take
$\mathcal{A}$ to be the algebra of linear transformations on
$\Sigma_0(k)$ or $\Sigma(k)$ viewed as vector spaces, we can take $x$
to be the linear transformation corresponding to multiplication by the
indeterminant $t$, and then $p(x)$ is the linear transformation
corresponding to multiplication by $p(t)$.

\section{Some other sums}
\label{some other sums}
\setcounter{equation}{0}

	Consider the vector space over ${\bf C}$ of doubly-infinite
sequences $\{a_j\}_{j=-\infty}^\infty$ of complex numbers such that
$a_j = 0$ when $|j|$ is sufficiently large, where the vector space
operations are termwise addition and scalar multiplication, as usual.
If $\{a_j\}_{j=-\infty}^\infty$, $\{b_l\}_{l=-\infty}^\infty$ are two
such sequences, then we can define the Cauchy product to be the
sequence $\{c_n\}_{n=-\infty}^\infty$ of complex numbers with
\begin{equation}
	c_n = \sum_{j=-\infty}^\infty a_j \, b_{n-j}.
\end{equation}
For each integer $n$, this sum has only finitely many nonzero, and so
it makes sense.  One can also check that $c_n = 0$ when $|n|$ is
sufficiently large, so that $\{c_n\}_{n=-\infty}^\infty$ lies in the
vector space under consideration.  In this way we get a commutative
algebra over the complex numbers.

	If $\{a_j\}_{j=-\infty}^\infty$ is a sequence of this type,
then we can define a function associated to it on the non-zero complex
numbers by
\begin{equation}
	\sum_{j=-\infty}^\infty a_j \, z^j.
\end{equation}
Again the sum makes sense because all but finitely many of the terms
are equal to $0$.  Thus we get a mapping from our algebra to functions
on ${\bf C} \backslash \{0\}$, and it is easy to see that this is an
algebra homomorphism, i.e., this mapping is linear and it takes Cauchy
products of sequences to ordinary products of functions on ${\bf C}
\backslash \{0\}$.

	Now consider doubly-infinite sequences of complex numbers
$\{a_j\}_{j=1}^\infty$ such that
\begin{equation}
	\sum_{j=-\infty}^\infty |a_j| < \infty,
\end{equation}
which is the same as saying that the partial sums $\sum_{j=-n}^n
|a_j|$ are bounded.  In other words the series corresponding to the
sequence should converge absolutely.  As before the space of these
sequences is a vector space over the complex numbers.  We can define a
norm on this space by saying that the norm of such a sequence is equal
to the sum of the moduli of its terms.  In this way our vector space
becomes a Banach space over the complex numbers, which is to say that
it is a normed vector space which is complete with respect to the
metric associated to the norm.

	If $\{a_j\}_{j=-\infty}^\infty$, $\{b_l\}_{l=-\infty}^\infty$
are two sequences of this type, with absolutely summable terms, then
we can define the Cauchy product $\{c_n\}_{n=-\infty}^\infty$ in the
same manner as before.  The series used to define $c_n$ converges
absolutely for each $n \in {\bf Z}$.  Moreover, one can show that
$\{c_n\}_{n=-\infty}^\infty$ also lies in our space, so that the
Cauchy product makes our vector space a commutative algebra over the
complex numbers.  The norm of the Cauchy product is less than or equal
to the product of the norms of the sequences used in the product, and
thus we have a commutative Banach algebra.

	If $\{a_j\}_{j=-\infty}^\infty$ is one of these
doubly-infinite sequences of complex numbers whose terms are
absolutely summable, then we can associate to it the function on the
unit circle given by
\begin{equation}
	\sum_{j=0}^\infty a_j \, z^j 
		+ \sum_{j=-\infty}^{-1} a_j \, \overline{z}^{-j}.
\end{equation}
This is equivalent to the earlier formula, because $z^{-1} =
\overline{z}$ when $|z| = 1$, and we need to restrict our attention to
complex numbers $z$ with $|z| = 1$ to ensure that the series
converges.  The absolute summability of the coefficients $a_j$ implies
that the partial sums of this series converge uniformly on the unit
circle as in the Weierstrass $M$-test, and hence the sum defines a
continuous function on the unit circle.  The mapping from sequences to
functions is again an algebra homomorphism, which is to say that it is
linear and takes Cauchy products of sequences to products of functions
on the unit circle.

\appendix

\addtocontents{toc}{{\protect\bigskip\protect\hbox{\protect\bf\protect\large Appendices}}}

\section{Metric spaces}
\label{metric spaces}
\setcounter{equation}{0}

\renewcommand{\thetheorem}{A.\arabic{equation}}
\renewcommand{\theequation}{A.\arabic{equation}}

	By a \emph{metric space} we mean a nonempty set $M$ together
with a real-valued function $d(x, y)$ defined for $x, y \in M$ such
that $d(x, y) \ge 0$ for all $x, y \in M$, $d(x, y) = 0$ if and only
if $x = y$, $d(x, y) = d(y, x)$ for all $x, y \in M$, and
\begin{equation}
	d(x, w) \le d(x, y) + d(y, w)
\end{equation}
for all $x, y, w \in {\bf M}$.  This last property is called the
\emph{triangle inequality} for the metric $d(x, y)$ on $M$.  If $(M,
d(x, y))$ is a metric space and $E$ is a nonempty subset of $M$, then
we can view $E$ as a metric space in its own right, using the
restriction of the metric $d(x, y)$ to $x, y \in M$.

	Suppose that $(M, d(x, y))$ is a metric space and that
$\{x_j\}_{j=1}^\infty$ is a sequence of points in $M$.  We say that
$\{x_j\}_{j=1}^\infty$ is a \emph{Cauchy sequence} in $M$ if for every
$\epsilon > 0$ there is a a positive integer $L$ such that
\begin{equation}
	d(x_j, x_l) < \epsilon
\end{equation}
for all $j, l \ge L$.  We say that $\{x_j\}_{j=1}^\infty$
\emph{converges} to a point $x \in M$ if for every $\epsilon > 0$
there is a positive integer $L$ such that
\begin{equation}
	d(x_j, x) < \epsilon
\end{equation}
for all $j \ge L$.  One can check that the limit of a convergent
sequence is unique, which is to say that if a sequence
$\{x_j\}_{j=1}^\infty$ of points in $M$ xonverges to $x \in M$ and to
$x' \in M$, then $x' = x$.  If $\{x_j\}_{j=1}^\infty$ is a sequence of
points in $M$ which converges to a point $x \in M$, then we write
\begin{equation}
	\lim_{j \to \infty} x_j = x.
\end{equation}

	It is easy to see that a convergent sequence in a metric space
is a Cauchy sequence.  Roughly speaking, the property of being a
Cauchy sequence captures the information of convergence without
having a limit.  A metric space is said to be \emph{complete} if every
Cauchy sequence in the space has a limit.

	Let $(M, d(x, y))$ be a metric space.  Suppose that
$\{x_j\}_{j=1}^\infty$, $\{y_j\}_{j=1}^\infty$ are sequences
in $M$ such that
\begin{equation}
	\lim_{j \to \infty} d(x_j, y_j) = 0
\end{equation}
as a sequence of real numbers.  One can show that if
$\{x_j\}_{j=1}^\infty$ is a Cauchy sequence, then the same is true of
$\{y_j\}_{j=1}^\infty$.  Similarly, if $\{x_j\}_{j=1}^\infty$
converges in $M$, then $\{y_j\}_{j=1}^\infty$ converges too, and to
the same point.

	A subset $E$ of $M$ is said to be \emph{dense} in $M$ if for
every point $x \in M$ there is a sequence $\{x_j\}_{j=1}^\infty$ of
points in $E$ which converges to $x$.  This is equivalent to saying
that for each $x \in M$ and each positive real number $\epsilon$ there
is a point $y \in E$ such that $d(x, y) < \epsilon$.  Indeed, when
this condition holds, one can simply choose $x_j \in E$ for each
positive integer $j$ such that $d(x_j, x) < 1/j$.

	If $(M, d(x, y))$ and $(N, \rho(u, v))$ are metric spaces,
then a mapping $\phi$ from $M$ into $N$ is said to be an
\emph{isometry} if
\begin{equation}
	\rho(\phi(x), \phi(y)) = d(x, y)
\end{equation}
for all $x, y \in M$.  Any metric space $(M, d(x, y))$ has a
\emph{completion} in the sense that there is a complete metric space
$(N, \rho(u, v))$ and an isometry $\phi : M \to N$ such that
\begin{equation}
	\phi(M) = \{\phi(x) : x \in M\}
\end{equation}
is dense in $N$.  The completion is unique up to isomorphism, in the
sense that if $(N', \rho'(z, w))$ is another complete metric space and
$\phi'$ is an isometric embedding of $M$ into $N'$ whose image is
dense in $N'$, then there is an isometry $\psi$ from $N$ onto $N'$
such that
\begin{equation}
	\psi(\phi(x)) = \phi'(x)
\end{equation}
for all $x \in M$.  Basically one defines $\psi$ first on the image of
$M$ in $N$ by this condition, and then shows that $\psi$ can be
extended to an isometry from $N$ onto $N'$.

	If $(M, d(x, y))$ is a metric space, then $d(x, y)$
is said to be an \emph{ultrametric} if
\begin{equation}
	d(x, w) \le \max (d(x, y), d(y, w))
\end{equation}
for all $x, y , w \in M$, which is a stronger version of the triangle
inequality.  For instance, suppose that $F_1, F_2, \ldots$, is a
sequence of finite sets, each with at least $2$ elements, and let $M$
be the set of sequences $x = \{x_l\}_{l=1}^\infty$ such that $x_l \in
F_l$ for each positive integer $l$.  Let $\{\rho_l\}_{l=1}^\infty$ be
a sequence of positive real numbers which is strictly decreasing and
converges to $0$.  If $x = \{x_l\}_{l=1}^\infty$, $y =
\{y_l\}_{l=1}^\infty$ are elements of $M$, then put $d(x, y) = 0$ if
$x = y$ and otherwise put $d(x, y) = \rho_n$, where $n$ is the
smallest positive integer such that $x_n \ne y_n$.  This defines an
ultrametric on $M$.

\section{Compactness}
\label{compactness}
\setcounter{equation}{0}

\renewcommand{\thetheorem}{B.\arabic{equation}}
\renewcommand{\theequation}{B.\arabic{equation}}

	Let $A$ be a set, and let $\{x_j\}_{j=1}^\infty$ be a sequence
of points in $A$.  By a \emph{subsequence} of $\{x_j\}_{j=1}^\infty$
we mean a sequence of the form $\{x_{j_l}\}_{l=1}^\infty$, where
$\{j_l\}_{l=1}^\infty$ is a strictly increasing sequence of positive
integers.  In other words, we basically restrict the original sequence
to an infinite subset of integer indices, arranged in increasing
order.  

	A sequence is automatically considered to be a subsequence of
itself, and a subsequence of a subsequence of a sequence
$\{x_j\}_{j=1}^\infty$ is also a subsequence of
$\{x_j\}_{j=1}^\infty$.  Observe that if $\{x_j\}_{j=1}^\infty$ is a
sequence of points in $A$ and if $A$ is the union of finitely many
subsets $A_1, \ldots, A_n$, then there is a $q$, $1 \le q \le n$, and
a subsequence $\{x_{j_l}\}_{l=1}^\infty$ of $\{x_j\}_{j=1}^\infty$
such that $x_{j_l} \in A_q$ for all positive integers $l$.

	Let $(M, d(x, y))$ be a metric space.  If
$\{x_j\}_{j=1}^\infty$ is a sequence of points in $M$ which converges
to a point $x \in M$, then every subsequence
$\{x_{j_l}\}_{l=1}^\infty$ of $\{x_j\}_{j=1}^\infty$ converges to $x$.
If $\{x_j\}_{j=1}^\infty$ is a Cauchy sequence in $M$, then every
subsequence $\{x_{j_l}\}_{l=1}^\infty$ of $\{x_j\}_{j=1}^\infty$ is
also a Cauchy sequence.

	A subset $E$ of $M$ is said to be \emph{closed} if for every
sequence $\{x_j\}_{j=1}^\infty$ of points in $E$ which converges to a
point $x \in M$ we have that $x \in M$.  This is equivalent to the
requirement that if $x$ is an element of $M$ such that for each
$\epsilon > 0$ there is a $y \in M$ with $d(x, y) < \epsilon$, then $y
\in E$.  The empty set and $M$ itself are automatically closed subsets
of $M$.

	If $A$ is a subset of $M$ and $p$ is a point in $M$, then that
$p$ is a \emph{limit point} of $A$ if for each $r > 0$ there exists a
point $x \in A$ such that $d(p, x) < r$ and $x \ne p$.  This is
equivalent to saying that for each $r > 0$ there are infinitely many
elements of $A$ whose distance to $p$ is less than $r$.  Equivalently,
$p$ is a limit point of $A$ if there is a sequence of points in $A$,
none of which are equal to $p$, and which converges to $p$.  At any
rate, a finite subset of $M$ has no limit points.  A subset of $M$ is
closed if and only if it contains all of its limit points.

	For each $p \in M$ and nonnegative real number $r$, the set
\begin{equation}
	\{x \in M : d(x, p) \le r\}
\end{equation}
is a closed subset of $M$.  One can show this using the triangle
inequality.  The union of finitely many closed subsets of $M$ is a
closed subset of $M$.  The intersection of any family of closed
subsets of $M$ is a closed subset of $M$.

	Let $A$ be an arbitrary subset of $M$.  The \emph{closure} of
$A$ is denoted $\overline{A}$ and is defined to be the set of points
$x \in M$ for which there is a sequence $\{x_j\}_{j=1}^\infty$ of
points in $A$ which converges to $x$.  If $x \in A$, then one can take
$x_j = x$ for all $j$, and thus we have that $A \subseteq
\overline{A}$.  Equivalently, a point $x \in M$ lies in the closure of
$A$ if for every $\epsilon > 0$ there is a point $y \in A$ such that
$d(x, y) < \epsilon$.  This is also the same as saying that the
closure of $A$ is equal to the union of $A$ and the set of limit
points of $A$.

	By definition, the closure of $A$ is equal to $A$ if and only
if $A$ is a closed set.  In general one can check that if $A$ is any
subset of $M$, then the closure of $A$ is a closed subset of $M$.  For
if $x$ is a point in $M$ which can be approximated by elements of
$\overline{A}$, then $x$ can be approximated by elements of $A$, and
therefore lies in $\overline{A}$.

	A subset $K$ of $M$ is said to be \emph{compact} if for every
sequence $\{x_j\}_{j=1}^\infty$ of points in $K$ there is a
subsequence $\{x_{j_l}\}_{l=1}^\infty$ of $\{x_j\}_{j=1}^\infty$ which
converges to a point in $K$.  A subset $K$ of $M$ has the \emph{limit
point property} if every infinite subset $E$ of $K$ has a limit point
which is an element of $K$.  It is not too difficult to show that
compactness is equivalent to the limit point property.  Namely, if $K$
is compact and $E$ is an infinite subset of $K$, one can choose a
sequence of points in $E$ in which no point occurs more than once, and
any subsequential limit of this sequence is a limit point of $E$.
Conversely, if $K$ has the limit point property and
$\{x_j\}_{j=1}^\infty$ is a sequence of points in $K$, then either
there is a subsequence of $\{x_j\}_{j=1}^\infty$ which is constant and
hence convergent, or the set of $x_j$'s is infinite and one can check
that a limit point of this set is also the limit of a subsequence of
$\{x_j\}_{j=1}^\infty$.

	A compact subset $K$ of $M$ is closed.  For suppose that $x$
is an element of $M$ and $\{x_j\}_{j=1}^\infty$ is a sequence of
points in $K$ which converges to $x$.  By compactness there is a
subsequence $\{x_{j_l}\}_{l=1}^\infty$ of $\{x_j\}_{j=1}^\infty$
which converges to a point in $K$.  However, this subsequence
converges to $x$, since the whole sequence converges to $x$,
and it follows that $x \in K$.

	If $K_1, \ldots, K_n$ are compact subsets of $M$, then the
union
\begin{equation}
	K = K_1 \cup K_2 \cup \cdots \cup K_n
\end{equation}
is also a compact subset of $M$.  Indeed, if $\{x_j\}_{j=1}^\infty$ is
a sequence of points in $K$, then there is a $q$, $1 \le q \le n$, and
a subsequence $\{x_{j_l}\}_{l=1}^\infty$ such that $x_{j_l} \in K_q$
for all $l$.  The compactness of $K_q$ then implies that there is a
subsequence of $\{x_{j_l}\}_{l=1}^\infty$ which converges to a point
in $K_q$.  This subsubsequence of $\{x_j\}_{j=1}^\infty$ is also a
subsequence of $\{x_j\}_{j=1}^\infty$, and it converges to a point in
$K$, as desired.

	Suppose that $E$, $K$ are subsets of $M$, with $E \subseteq
K$, $E$ closed, and $K$ compact.  Under these conditions $E$ is also
compact.  For suppose that $\{x_j\}_{j=1}^\infty$ is a sequence of
points in $E$.  Because $K$ is compact, a subsequence of
$\{x_j\}_{j=1}^\infty$ converges to a point in $K$.  Because $E$ is
closed, the limit of this subsequence lies in $E$.

	A subset $E$ of $M$ is said to be \emph{bounded} if there is a
positive real number $r$ such that $d(x, y) \le r$ for all $x, y \in
E$.  Equivalently, $E$ is bounded if there exists $p \in M$ and $t >
0$ such that $d(p, x) \le t$ for all $x \in E$.  This is also
equivalent to the condition that for each $p \in M$ there is a $t(p) >
0$ so that $d(p, x) \le t(p)$ for all $x \in M$.  Observe that the
union of finitely many bounded sets is bounded.

	If $E$ is a nonempty bounded subset of $M$, then the
\emph{diameter} of $E$ is denoted $\diam E$ and defined to be the
supremum of the numbers $d(x, y)$ for $x, y \in E$.  The closure of a
bounded set is also bounded, and has the same diameter, assuming that
it is not empty.  The diameter of a union of two bounded subsets of
$M$ is less than or equal to the sum of the diameters of the two
subsets if the two subsets have a point in common, and if the metric
is an ultrametric, then the diameter of the union is less than or
equal to the maximum of the diameters of the two subsets.

	A compact subset $K$ of $M$ is bounded.  Indeed, let $p$ be
any element of $M$.  If $K$ is not bounded, then for each positive
integer $n$ there is a point $x_n \in K$ such that $d(p, x_n) \ge n$.
In this case $\{x_n\}_{n=1}^\infty$ is a sequence of points in $K$ for
which there is no convergent subsequence, contradicting the assumption
that $K$ is compact.

	A subset $E$ of $M$ is said to be \emph{totally bounded} if
for each $\epsilon > 0$ there exist finitely many points $p_1, \ldots,
p_l \in M$ so that for each $x \in K$ we have $d(x, p_j) \le \epsilon$
for at least one $j$.  Equivalently, $E$ is totally bounded if it can
be expressed as the union of finitely many subsets of arbitrarily
small diameter.  The closure of a totally bounded subset of $M$ is
also totally bounded.

	Compact subsets of $M$ are totally bounded.  Indeed, assume
for the sake of a contradiction that $K$ is a compact subset of $M$
which is not totally bounded.  In this case there is an $\epsilon > 0$
and a sequence of points $\{x_l\}_{l=1}^\infty$ in $K$ such that
$d(x_l, x_m) \ge \epsilon$ when $l < m$.  Clearly no subsequence of
$\{x_l\}_{l=1}^\infty$ can converge.  In fact, no subsequence of
$\{x_l\}_{l=1}^\infty$ can be a Cauchy sequence.

	If $E$ is a subset of $M$ which is totally bounded, and if
$\{x_l\}_{l=1}^\infty$ is any sequence of points in $E$, then for each
$\epsilon > 0$ there is a subsequence of $\{x_l\}_{l=1}^\infty$
contained in a subset of $E$ of diameter less than $\epsilon$.  This
is easy to see by expressing $E$ as a union of finitely many subsets
each with diameter less than $\epsilon$.  One can go a bit further and
say that if $E$ is a totally bounded subset of $M$ and
$\{x_l\}_{l=1}^\infty$ is a sequence of points in $E$, then
$\{x_l\}_{l=1}^\infty$ has a subsequence which is a Cauchy sequence.
This uses a Cantor diagonalization argument.

	In short, a subset $E$ of $M$ is totally bounded if and only
if every sequence of points in $E$ has a subsequence which is a Cauchy
sequence.  If $M$ is a complete metric space, then a subset of $M$ is
compact if and only if it is closed and totally bounded, by the
preceding observations.  In any metric space, a Cauchy sequence with a
convergent subsequence converges to the same limit as the subsequence,
and thus a Cauchy sequence contained in a compact set converges.  It
follows that a metric space $M$ is compact as a subset of itself if
and only if it is complete and totally bounded.

	An interesting class of examples is provided by the spaces
mentioned at the end of Appendix \ref{metric spaces}, consisting of
sequences with the $l$th term in a fixed finite set $F_l$ for each
$l$.  For these spaces one can verify completeness, total boundedness,
and compactness quite concretely.  A sequence of elements in one of
these spaces is a sequence of sequence, and convergence basically
amounts to convergence of the individual terms in the $F_l$'s.

\section{Topological spaces}
\label{topological spaces}
\setcounter{equation}{0}

\renewcommand{\thetheorem}{C.\arabic{equation}}
\renewcommand{\theequation}{C.\arabic{equation}}

	Let $X$ be a set, and let $\tau$ be a collection of subsets of
$X$, called the \emph{open} subsets of $X$.  We say that $\tau$
defines a \emph{topology} on $X$, so that $X$ becomes a
\emph{topological space}, if the empty set $\emptyset$ and $X$ itself
are open subsets of $X$, if the intersection of finitely many open
subsets of $X$ is again an open subset of $X$, and if the union of any
family of open subsets of $X$ is an open subset of $X$.  The condition
about unions is equivalent to saying that if $W$ is a subset of $X$
and if for each $p \in W$ there is an open subset $U$ of $X$ such that
$p \in U$ and $U \subseteq W$, then $W$ is an open subset of $X$.

	Let $(M, d(x, y))$ be a metric space.  If $p \in M$ and $r >
0$, then the open ball with center $p$ and radius $r$ is denoted $B(p,
r)$ and defined to be the set of $x \in M$ such that $d(p, x) < r$.
The closed ball with center $p$ and radius $r$ is denoted
$\overline{B}(p, r)$ and defined to be the set of $x \in M$ such that
$d(p, x) \le r$.  A subset $U$ of $M$ is said to be open if for each
$p \in U$ there is an $r > 0$ such that $B(p, r) \subseteq U$.  It is
easy to check that this defines a topology on $M$ and that every
open ball in $M$ is an open subset of $M$.

	If $(X, \tau)$ is a topological space, then a subset $E$ of
$X$ is said to be \emph{closed} if the complement $X \backslash E$ of
$E$ in $X$, consisting of the points in $X$ which do not lie in $E$,
is an open subset of $X$.  It follows that the empty set and $X$
itself are closed subsets of $X$, that the union of finitely many
closed subsets of $X$ is a closed subset of $X$, and that the
intersection of any family of closed subsets of $X$ is a closed subset
of $X$.  In a metric space a subset is closed in the sense defined in
Appendix \ref{compactness} if and only if it is closed in the sense
that its complement is open.

	If $A$ is an arbitrary subset of $X$, then the \emph{closure}
of $A$ is denoted $\overline{A}$ and defined to be the set of points
$p \in X$ such that for each open subset $U$ of $X$ with $p \in U$ we
have that the intersection of $A$ and $U$ is nonempty.  Thus the
closure of $A$ contains $A$ automatically, and in a metric space this
definition of the closure is equivalent to the one in Appendix
\ref{compactness}.  One can check that $A = \overline{A}$ if and only
if $A$ is a closed subset of $X$, and that $\overline{A}$ always is a
closed subset of $X$.

	If $X$, $Y$ are topological spaces and $f$ is a mapping from
$X$ to $Y$, then we say that $f$ is \emph{continuous} at a point $p
\in X$ if for each open subset $W$ of $Y$ such that $f(p) \in W$ there
is an open subset $U$ of $X$ such that $p \in U$ and
\begin{equation}
	U \subseteq f^{-1}(W).
\end{equation}
Recall that $f^{-1}(W)$ is by definition the set of points $x \in X$
such that $f(x) \in W$.  If $f : X \to Y$ is continuous at every point
in $X$, then we say that $f$ is a continuous mapping from $X$ to $Y$.
This is equivalent to saying that $f^{-1}(V)$ is an open subset of $X$
for every open subset $V$ of $Y$, or that $f^{-1}(E)$ is a closed
subset of $X$ for every closed subset $E$ of $Y$.

	Let $X$ and $Y$ be sets, and let $f$ be a mapping from $X$ to
$Y$.  We say that $f$ is injective or one-to-one if for each pair of
points $x_1, x_2 \in X$ with $x_1 \ne x_2$ we have that $f(x_1) \ne
f(x_2)$.  We say that $f$ maps $X$ onto $Y$ if for each $y \in Y$
there is an $x \in X$ such that $f(x) = y$.  Thus $f$ is a one-to-one
mapping of $X$ onto $Y$ if and only if there is an inverse mapping $h$
from $Y$ to $X$ such that
\begin{equation}
	h(f(x)) = x
\end{equation}
for all $x \in X$ and
\begin{equation}
	f(h(y)) = y
\end{equation}
for all $y \in Y$.  The inverse mapping is unique when it exists, and
is denoted $f^{-1}$.

	Let $X$, $Y$ be topological spaces, and let $f$ be a
one-to-one mapping from $X$ onto $Y$.  We say that $f$ is a
\emph{homeomorphism} if $f$ is a continuous mapping from $X$ to $Y$
and $f^{-1}$ is a continuous mapping from $Y$ to $X$.  Equivalently,
$f$ is a homeomorphism if $f$ sends open subsets of $X$ to open
subsets of $Y$ and $f^{-1}$ sends open subsets of $Y$ to open subsets
of $X$.  One could just as well use closed subsets instead of open
subsets here.

	A subset $K$ of a topological space $X$ is said to be
\emph{compact} if every open covering of $K$ admits a finite
subcovering.  Recall that an \emph{open covering} of a subset $K$ of
$X$ is a family $\{U_\alpha\}_{\alpha \in A}$ of open subsets of $X$
whose union contains $K$ as a subset.  Thus $K$ is compact if for each
open covering $\{U_\alpha\}_{\alpha \in A}$ of $K$ there is a finite
subset $A_1$ of $A$ such that $K$ is contained in the union of the
$U_\alpha$'s with $\alpha \in A_1$.  Finite subsets of $X$ are
automatically compact.

	If $X$, $Y$ are topological spaces, $f$ is a continuous
mapping from $X$ to $Y$, and $K$ is a compact subset of $X$, then
$f(K)$ is a compact subset of $Y$, where $f(K)$ is by definition the
set of points in $Y$ of the form $f(x)$ for some $x \in K$.  Indeed,
suppose that $\{V_\alpha\}_{\alpha \in A}$ is an arbitrary open
covering of $f(K)$ in $Y$.  Then $\{f^{-1}(V_\alpha)\}_{\alpha \in A}$
is an open covering of $K$ in $X$.  Because $K$ is compact, there is a
finite subset $A_1$ of $A$ such that $K$ is contained in the union of
$f^{-1}(V_\alpha)$, $\alpha \in A_1$.  This implies that $f(K)$ is
contained in the union of $V_\alpha$, $\alpha \in A_1$.

	A topological space $X$ is said to satisfy the first axiom of
separation if subsets of $X$ with exactly one element are closed
subsets.  This is equivalent to saying that if $x, y \in X$ and $x \ne
y$, then there is an open subset $U$ of $X$ such that $x \in U$ and $y
\not\in U$.  If $X$ satisfies the first axiom of separation, then
every finite subset of $X$ is a closed subset of $X$.  A topological
space $X$ is said to satisfy the second axiom of separation if for
every $x, y \in X$ with $x \ne y$ there are open subsets $U$, $V$ of
$X$ such that $x \in U$, $y \in V$, and $U \cap V = \emptyset$.  One
also calls $X$ a Hausdorff topological space in this case, and of
course the second axiom of separation implies the first axiom of
separation.

	Let $X$ be a Hausdorff topological space, let $K$ be a compact
subset of $X$, and fix a point $x \in X \backslash K$.  For each $y
\in K$ there are open subsets $U(y)$, $V(y)$ of $X$ such that $x \in
U(y)$, $y \in V(y)$, and $U(y)$, $V(y)$ are disjoint.  Because $K$ is
compact there are finitely many points $y_1, \ldots, y_n \in K$ such
that
\begin{equation}
	K \subseteq V(y_1) \cup \cdots \cup V(y_n).
\end{equation}
Thus
\begin{equation}
	U(y_1) \cap \cdots U(y_n)
\end{equation}
is an open subset of $X$ which contains $x$ and is contained in the
complement of $K$.  It follows that $K$ is a closed subset of $X$.

	In order for compact subsets of a topological space $X$ to be
closed, it is obviously necessary for $X$ to satisfy the first axiom
of separation, since finite subsets of $X$ are compact.  The first
axiom of separation is not sufficient, as one can show by examples.
For this reason the term ``quasicompact'' is sometimes used for
compact subsets in topological spaces which are not necessarily
Hausdorff.

	For instance, let $X$ be the real line ${\bf R}$ together with
an extra point $0'$.  Let us write ${\bf R}'$ for the set obtained
from ${\bf R}$ by removing $0$ and adding $0'$.  One can define a
topology on $X$ so that ${\bf R}$ and ${\bf R}'$ are open subsets of
$X$ which are each homeomorphic to the real line with the standard
topology.  The resulting space $X$ satisfies the first axiom of
separation but not the second one, because if $U$, $V$ are open
subsets of $X$ such that $0 \in U$ and $0' \in V$, then $U$ and $V$
have elements in common, namely nonzero real numbers close to $0$.
Closed and bounded intervals in the real line are compact, and they
give rise to compact subsets of ${\bf R}$ and ${\bf R}'$ which are
compact subsets of $X$, but they may not be closed because of
containing $0$ and not $0'$ or vice-versa.

	In any topological space $X$, the union of two compact sets
$K_1$, $K_2$.  For suppose that $\{U_\alpha\}_{\alpha \in A}$ is an
open covering of $K_1 \cup K_2$.  Then this is also an open covering
of $K_1$, $K_2$ individually, and thus there are finite subsets $A_1$,
$A_2$ of $A$ such that $K_1$ is contained in the union of the
$U_\alpha$'s with $\alpha \in A_1$ and $K_2$ is contained in the union
of the $U_\alpha$'s with $\alpha \in A_2$.  Therefore $A_1 \cup A_2$
is a finite subset of $A$ such that $K_1 \cup K_2$ is contained in the
union of the $U_\alpha$'s with $\alpha \in A_1 \cup A_2$.

	In any topological space $X$, if $F$ is a closed subset of $X$
and $K$ is a compact subset of $X$, then the intersection of $F$ and
$K$ is compact.  Indeed, let an arbitrary open covering of $F \cap K$
be given.  Since $F$ is closed, $X \backslash F$ is open, and we can
add this open set to the open covering of $F \cap K$ to get an open
covering of $K$.  The compactness of $K$ implies that finitely many
open subsets in the original open covering of $F \cap K$ together
perhaps with $X \backslash F$ covers $K$.  Hence these finitely many
open subsets from the original covering of $F \cap K$ covers $F \cap
K$.

	If $(M, d(x, y))$ is a metric space, then a subset $K$ of $M$
is compact in the sense of open coverings if and only if it is compact
in the sense of Appendix \ref{compactness}.  To see this, assume first
that $K$ is compact in the sense of open coverings, and let $E$ be any
infinite subset of $K$.  If we assume for the sake of a contradiction
that $E$ does not have a limit point in $K$, then for each $p \in K$
there is then a positive real number $r(p)$ such that $E \cap B(p,
r(p))$ is finite.  Since $K$ is compact in the sense of open
coverings, $K$ is contained in the union of $B(p, r(p))$ for finitely
many $p \in K$, which implies that $E$ is finite.  Notice also that
compactness in the sense of open coverings immediately implies the
property of being totally bounded.

	Conversely, suppose that $K$ is a totally bounded subset of
$M$ and that every Cauchy sequence in $K$ converges to a point in $K$.
Suppose too that there is an open covering of $K$ which does not admit
any finite subcovering.  Using the total boundedness of $K$ one can
show that there is a sequence $K_1, K_2, \ldots$ of closed subsets of
$K$ such that $K_{n+1} \subseteq K_n$ for all $n$, the diameter of
$K_n$ is less than $1/n$ for all $n$, and no finite subcollection of
open sets from our covering of $K$ covers any $K_n$.  Using the
convergence of Cauchy sequences one can show that there is a point $p
\in K$ such that $p \in K_n$ for all $n$.  Hence $p$ is contained in
one of the open subsets in the covering of $K$, which then contains
$K_n$ for sufficiently large $n$, a contradiction.

\end{document}